\documentclass[twocolumn,english,APS,PRE,eqsecnum]{revtex4-2}
\usepackage[T1]{fontenc}
\usepackage[latin9]{inputenc}
\usepackage{color}
\usepackage{babel}
\usepackage{float}
\usepackage{fancybox}
\usepackage{calc}
\usepackage{bm}
\usepackage{varwidth}
\usepackage{amsmath}
\usepackage{amssymb}
\usepackage{graphicx}
\usepackage{geometry}
\usepackage{esint}
\usepackage[bookmarks=true,bookmarksnumbered=false,bookmarksopen=false,
 breaklinks=false,pdfborder={0 0 1},backref=false,colorlinks=false]
 {hyperref}

\makeatletter

\providecommand{\tabularnewline}{\\}
\newenvironment{cellvarwidth}[1][t]
    {\begin{varwidth}[#1]{\linewidth}}
    {\@finalstrut\@arstrutbox\end{varwidth}}

\makeatletter
\def\ps@pprintTitle{%
  \let\@oddhead\@empty
  \let\@evenhead\@empty
  \def\@oddfoot{\hfil\thepage\hfil}
  \let\@evenfoot\@oddfoot
}
\makeatother

\usepackage{xcolor}
\definecolor{BLUE}{RGB}{0,0,255}
\usepackage{array}

\makeatother

\begin{document}
\title{Non-periodic Fourier propagation algorithms for partial differential
equations}
\author{Channa Hatharasinghe$^{1}$, Run Yan Teh$^{1}$, Jesse van Rhijn$^{2}$,
Margaret D. Reid$^{1}$, Peter D. Drummond$^{1}$}
\address{$^{1}$Centre for Quantum Science and Technology Theory, Swinburne
University of Technology, Melbourne, Australia}
\address{$^{2}$University of Twente, Enschede, The Netherlands.}
\begin{abstract}
Spectral methods for partial differential equations (PDEs) with non-periodic
boundary conditions arising in computational physics often use polynomial
expansions on non-uniform grids. Here, we implement a Fourier method
that employs fast trigonometric expansions on a uniform grid with
non-periodic boundaries using fast discrete sine transforms (DST)
or/and discrete cosine transforms (DCT) to solve parabolic PDEs. We
implement this method in two ways: either using a Fourier spectral
derivative or a Fourier interaction picture. Both methods can treat
vector fields with a combination of Dirichlet and/or Neumann boundary
conditions in one or more space dimensions. As examples, we use them
to solve a variety of computational physics PDEs with analytical solutions,
including the Peregrine solitary wave solution. For the 1D heat equation
problem, our method with an interaction picture is accurate up to
machine precision. Soluble examples of stochastic partial differential
equation (SPDE) with non-periodic boundaries in one and two space
dimensions, with physics and interdisciplinary applications are also
treated. We compare the results obtained from these algorithms with
publicly available solvers that use polynomial spectral methods, and
study their relative performance and error scaling. Polynomial methods
with non-uniform spatial grids have lower spatial discretization errors
when the solutions change slowly in space, typically with large spatial
grids. For problems with rapid spatial variation, Fourier methods
can outperform polynomial expansions, owing to their smaller maximum
space interval, and are generally faster due to the computational
efficiency of discrete Fourier transform methods. We verified this
by making a complexity analysis in which we studied the total error
at the optimum combination of time and space steps for a given resource
use. This relative advantage was particularly apparent for stochastic
problems with delta-correlated spatial noise, in which there is no
asymptotically smooth regime. We also demonstrated applications to
stochastic problems in higher space dimensions, and to physically
relevant spatial noise having finite correlation lengths, which is
efficiently generated using Fourier transform filtering on a uniform
spatial grid. 
\end{abstract}
\maketitle

\section{Introduction}

Partial differential equations (PDEs) occur widely in physics and
related fields \citep{courant2024methods,evans2022partial}. Spectral
transforms are a common technique for solving them numerically \citep{Bernardi1997,boyd2001,Canuto1988,Hesthaven2007,Kopriva2009,Canuto2006Spectral,Grandclement2009}.
These convert a PDE into a series of smoothly changing basis functions,
which typically are either trigonometric \citep{Hesthaven2007} or
polynomial functions \citep{Johnson1996ChebyshevPI,Weideman1999,Grandclement2006}.
Fourier methods are widely used for problems with periodic boundary
conditions. Such methods are also used for stochastic partial differential
equations (SPDEs) \citep{WERNER1997,Tyrrell2005pseudospectral}. However,
PDEs and SPDEs in physics, as well as in stochastic bridges and quantum
foundation theory, can also have non-periodic boundary conditions
\citep{khoshnevisan2014analysis,drummond2017forward,drummond2021time,reid2022quantum}.

In this paper, we describe and implement two Fourier based spectral
methods for integrating parabolic and hyperbolic partial differential
equations (PDEs) and stochastic partial differential equations (SPDEs)
on non-periodic boundaries. These methods use efficient trigonometric
spectral transforms, namely discrete sine transform (DST) and discrete
cosine transforms (DCT) \citep{krylov1906approximate,Jain1976,Pueschel2003,Strang1999,Shao2008}.
There are some previous studies \citep{Kosloff1983A_non_periodic_fourier,Wise2021pseudospectral},
but our interaction picture algorithm has clear advantages over these.
Fourier-based interaction-picture strategies with periodic boundaries
have been useful in physics applications, such as the nonlinear Schrödinger
equation \citep{Weideman1986,Taha1984}. The advantage of extending
this to non-periodic cases is that these methods are accurate, fast
and well-suited to the use of equidistant collocation points. This
reduces errors when the solution has localized peaks or stochasticity.
We demonstrate our approach and compare it with finite differences,
and with polynomial spectral methods used in two public domain software
libraries.

To make a brief comparison of methods for boundary value problems
in physics, PDEs are often treated using a discretized finite difference
representation of the differential operator \citep{Langtangen2017}.
This method can have convergence and stability problems since the
highest-order differential terms have large eigenvalues, which require
small time steps \citep{tadmor2012review,Douglas1961-ua}. Generally,
as we show in the Appendix, this method is rather inaccurate and slow.
The second-order differential operator, the Laplacian ($\nabla^{2}$),
consists of eigenvalues that grow quadratically, while the fourth-order
differential operator, the biharmonic ($\nabla^{4}$), has eigenvalues
that grow as a quartic. To ensure accuracy, it is necessary to use
finer grids to capture the behavior of the solution. The time step
must then be kept very small to maintain stability.

The Crank-Nicolson finite difference method \citep{Crank1947} works
well in one spatial dimension and can reduce these instabilities.
For higher dimensions, the complexity of this algorith reduces the
method's efficiency \citep{Duffy2004,Sharma2024,Britz2003}. Alternating
direction implicit methods (ADI) \citep{Birkhoff1962} are therefore
used, which still require very small time steps \citep{Sun2004}.
To solve these problems, finite difference methods can be replaced
with spectral methods, which have proven to produce accurate results
in many cases \citep{Fletcher1984,IbarraVillalon2023,Hult2007,Balac2016}.
For computational physics problems with periodic boundary conditions,
Fourier spectral methods on a uniform grid are frequently utilized.
In deterministic cases with non-periodic PDEs, polynomial basis functions
of orthogonal polynomials are commonly used. Either method has been
employed in many areas of physics, for example, to solve the quantum
stochastic nonlinear Schrödinger equation \citep{drummond1993simulation,drummond1993quantum,corney2006many}
and the Gross-Pitaevskii equation (GPE) either using Fourier transforms
with periodic boundaries \citep{Javanainen2006} or polynomials with
non-periodic boundaries \citep{blakie2008numerical,Wang2011}.

In cases without periodicity, the most common pseudo-spectral methods
use polynomial expansions, typically with a non-uniform spatial grid
\citep{Fornberg1996,Orszag1972}. This has a substantial numerical
cost, and is less useful for solving stochastic partial differential
equations (SPDE) because the stochastic noise variance increases with
inverse step size, which makes it more difficult to check for convergence
in the spatial discretization. The non-uniformity problem is overcome
by using a uniform grid, as in the spectral techniques described here.
It is also simpler to compare numerical and reference solutions on
a uniform grid, making it useful in error analysis and deep learning
applications for PDEs \citep{liu2023spfno,Kovachki2024Neuraloperator,Liu2023OPNO},
as well as being more suitable for spatially correlated noise.

We show that this type of Fourier approach can accommodate non-periodic,
time-varying Dirichlet, Robin, and Neumann boundary conditions for
arbitrary spatial dimensions and field components. The boundary type
modifies the transformation used, and can affect an equation's error
stability properties, so we consider a combination of various boundary
conditions labelled by an index $\beta=1,\ldots4$:
\begin{enumerate}
\item Dirichlet-Dirichlet (D-D)
\item Dirichlet-Neumann (D-N)
\item Neumann-Dirichlet (N-D)
\item Neumann-Neumann (N-N).
\end{enumerate}
These boundaries can be combined with periodic boundaries in other
dimensions, which we label as $\beta=0$ if needed. This allows different
boundary conditions for each field component, dimension, and boundary.
This can have advantages over polynomial approaches with unequal collocation
spacings when spatial derivatives need to be computed, because Fourier
methods are often faster, due to the intrinsic $N\ln N$ speed of
FFT transforms \citep{cooley1965algorithm,Frigo1998,Frigo2005}, combined
with a diagonal propagation matrix in the interaction picture.

Unequal grid spacings in some polynomial methods lead to larger spatial
steps in the center of the interval, which causes higher spatial discretization
errors when there are localized peaks near the centre of an interval,
as found in many PDEs and SPDEs. The two alternative Fourier methods
we analyze in detail here for solving such equations with first-order
time derivatives and non-periodic boundaries are:
\begin{description}
\item [{FSD}] A Fourier spectral derivative method.
\item [{FIP}] A Fourier interaction picture method.
\end{description}
These are a special case of the Galerkin approximation method \citep{Galerkin1968}.
This equivalence is demonstrated in Appendix A and Appendix B for
solving a specific type of problem: the 1D heat equation with non-periodic
boundaries. This equivalence is not limited to 1D problems; it can
also be extended to higher dimensions provided suitable patch functions
are available.

We implement this approach in a public domain software toolbox \citep{Kiesewetter2016,Kiesewetter2023,Drummond2023}.
We also carry out comparison studies on practical examples to show
how the results differ from previous algorithms. While not exhaustive,
this shows where spectral methods have advantages and disadvantages.
Next, we compare this approach with two other publicly available
spectral PDE solvers, namely Dedalus \citep{Burns2020}, a public-domain
package that implements polynomial spectral methods on a non-uniform
grid using a sparse Tau method, and Pdepe \citep{Skeel1990}, a function
in the MATLAB language that uses a Galerkin method \citep{Galerkin1968}
to discretize a PDE on a uniform grid. In the tables in sections \ref{sec:Comparisons-between-uniform}
and \ref{sec:stochastic-heat}, these are abbreviated as ``Tau''
and ``Gal'' respectively.

As examples, we consider PDE cases with time-dependent Dirichlet and/or
Neumann boundary conditions. Our results show that trigonometric spectral
methods typically have lower spatial discretization errors compared
to polynomial methods for nonlinear problems with localized regions
of high spatial gradients, and can be faster. For smoothly varying
solutions the polynomial spectral methods can give lower spatial discretization
errors, and this issue is taken up in detail using complexity order
analysis. When using the FIP for linear problems, such as solving
the heat equation, the Fourier methods can achieve accuracy up to
machine precision. The advantages of higher accuracy and faster speed
is especially important with SPDEs. These have high spatial variation,
and require repeated integrations to give low sampling errors to obtain
the optimum complexity exponent \citep{teh2025complexity}.

The layout of the paper is as follows. Section \ref{sec:Fourier-methods-with}
discusses the use of Fourier methods for PDEs on a uniform grid with
non-periodic boundaries. In section \ref{sec:Inhomogeneous-patch-solutions},
we discuss how these methods can be applied to solve PDEs with inhomogeneous
boundaries. Section \ref{sec:Comparisons-between-uniform} gives the
relative error obtained with our methods, compared with those obtained
with public domain software using polynomial expansions for some test
problems with exact solutions. The complexity order is explained and
computed. In section \ref{sec:stochastic-heat}, we compute the solution
of a stochastic partial differential equation (SPDE), compared with
an analytic solution, using the Fourier and polynomial spectral methods.
Error and time scaling properties are analyzed for the 1D stochastic
heat equation, as well as 1D and 2D stochastic heat equations with
filtered noise. A conclusion and a summary of results are given in
Section \ref{sec:Conclusion}.

\section{Fourier methods with non-periodic boundaries\label{sec:Fourier-methods-with}}

\subsection{The partial differential equation}

We wish to solve an initial value problem for a $c$-component real
or complex vector field $\bm{u}=\left[u^{1},\ldots u^{c}\right]$
that satisfies a parabolic, possibly stochastic partial differential
equation (PDE) in $1+d$ dimensions where $\bm{r}=\left[t,x_{1},\ldots x_{d}\right]$,
with a finite interval in each dimension, of the form
\begin{equation}
\begin{split}\frac{\partial\boldsymbol{u}}{\partial t}=\bm{g}\left[\bm{r},\bm{u}\right]+\bm{\mathcal{L}}\left[\boldsymbol{u}\right]\end{split}
\label{eq:pde}
\end{equation}
where $\bm{\mathcal{L}}$ is a differential term defined below, and
\begin{equation}
\begin{split}\bm{g}\left[\bm{r},\bm{u}\right]=\mathbf{A}\left[\bm{r},\bm{u}\right]+\underline{\mathbf{B}}\left[\bm{r},\bm{u}\right]\cdot\mathbf{w}(\mathbf{r})\end{split}
.\label{eq:deriv_without_linear_term-1}
\end{equation}

Here, $\mathbf{A}\left[\bm{r},\bm{u}\right]$ is a linear or nonlinear
term which in general may also include derivatives. In the SPDE case
$\underline{\mathbf{B}}\left[\bm{r},\bm{u}\right]$ is a multiplicative
stochastic noise matrix, while $\mathbf{w}(\mathbf{r})$ is a Gaussian
stochastic noise field such that
\begin{equation}
\left\langle w_{i}\left(\bm{r}\right)w_{j}\left(\bm{r}'\right)\right\rangle =\delta_{ij}\delta\left(t-t'\right)\delta^{d}\left(\bm{x}-\bm{x}'\right).
\end{equation}
 Stratonovich calculus is assumed here, so that no further Ito corrections
are required, provided the multiplicative stochastic noise matrix
is evaluated at the midpoint \citep{DRUMMOND1991,WERNER1997}.

The quasi-linear operator $\bm{\mathcal{L}}$ may include several
differential terms. Any order derivative is possible , where:
\begin{equation}
\bm{\mathcal{L}}\left[\boldsymbol{u}\right]^{\nu}=\sum_{\mu,o,i}D_{o,i}^{\nu,\mu}\frac{\partial^{o}}{\partial x_{i}^{o}}u^{\mu}.\label{eq:linear-derivative}
\end{equation}
Here $o>0$, $i=1,\ldots d$ , $\mu,\nu=1,\ldots c$. The matrix term
D can be a function of position $\bm{r}$ and field vector $\bm{u}$
in the FSD case, but must be constant in the more restrictive FIP
case (see also \ref{subsec:Algorithm}). Combinations of the two algorithms
can be used, by treating some of the terms using the FIP and others
using FSD. This is useful if some derivatives have constant coefficients,
while others have variable coefficients that depend on $\bm{r},\bm{u}$,
and therefore are more suitable to the FSD method. Cross-derivatives
can be included between different space dimensions, but are not treated
here to simplify the notation.

\subsection{Boundary types\label{subsec:Boundary-types}}

The boundary conditions treated are functionals of the derivatives
and boundary values, defined such that:
\begin{equation}
\mathcal{B}\left[\bm{r},\boldsymbol{\nabla},\bm{u}\right]=0.
\end{equation}
 The boundary space $\Omega$ is the surface of a hyper-rectangle.
For dimension $i$, boundary conditions are applied at $X_{i}^{(a,b)}$,
which define the $2d$ boundaries of a subspace of $\mathbb{R}^{d}$,
each having $d-1$ dimensions, where $\bm{r}=\left(t,\bm{x}\right)$.
The subspace is therefore defined by: 
\begin{equation}
X_{i}^{\left(a\right)}\le x_{i}\le X_{i}^{\left(b\right)}.
\end{equation}

The boundary conditions considered are of three types. These may be
applied independently for each field component and at each upper and
lower boundary. The fundamental boundary types treated are:
\begin{itemize}
\item Dirichlet (specified value):
\begin{equation}
u^{\nu}\left(\bm{r}_{a,b}^{(i)}\right)=U_{a,b}^{\nu i}\left(\bm{r}_{a,b}^{(i)}\right).
\end{equation}
\item Neumann (specified derivative):
\begin{equation}
\partial^{i}u^{\nu}\left(\bm{r}_{a,b}^{(i)}\right)=N_{a,b}^{\nu i}\left(\bm{r}_{a,b}^{(i)}\right).
\end{equation}
\item Periodic (periodic value and periodic derivative):
\begin{align}
u^{\nu}\left(\bm{r}_{a}^{(i)}\right) & =u^{\nu}\left(\bm{r}_{b}^{(i)}\right)\nonumber \\
\partial^{i}u^{\nu}\left(\bm{r}_{a}^{(i)}\right) & =\partial^{i}u^{\nu}\left(\bm{r}_{b}^{(i)}\right).
\end{align}
\end{itemize}
There are $5$ combinations possible for boundaries in each field
component and space dimension, since Dirichlet and Neumann boundaries
can be combined in four ways in a given component and dimension. These
cannot be used with periodic boundaries for the same field component
and dimension. This gives a total of $5^{cd}$ distinct boundary types,
each with boundary values that may be time-dependent. We label these
boundary types $\beta_{i}^{\mu}$ for component $\mu$ and dimension
$i$. Since periodic examples are treated elsewhere \citep{drummond1993simulation,WERNER1997,Tyrrell2005pseudospectral},
we focus on the four non-periodic cases, $\beta=1,\ldots4$.

More generally the boundary value vector functions $\bm{U}$ and $\bm{N}$
can be nonlinear functions of the fields, which requires modifications
to the methods described here. The definitions given above allow the
boundary values to have both time and space dependence. As a result
these may have stochastic or time-dependent values.

The boundary values themselves can be the solutions of subsidiary
differential or stochastic differential equations, which must then
be treated in parallel with the main partial differential equation
defined in Eq. (\ref{eq:pde}). For simplicity we do not include the
additional solutions that would be required in such more complex cases.

A general PDE problem in the form treated here can have inhomogeneous
boundary conditions. As the first step, the problem is converted into
a homogenous boundary value problem by reformulating the problem.
Similar techniques are used for polynomial spectral methods \citep{boyd2001}.
We use a change of variables $\mathbf{u}=\mathbf{v}+\bm{h}$, where $\mathbf{v}$ is
a patch function \citep{neidinger2019multivariate} chosen to satisfy
the inhomogeneous boundary conditions of the problem over a short
interval, so that we now must solve a new PDE using transforms of
the homogeneous part $\bm{h}$ 
\begin{equation}
\frac{\partial\boldsymbol{u}}{\partial t}=\mathbf{g}\left[\bm{r},\bm{u}\right]+\bm{\mathcal{L}}\left(\mathbf{\boldsymbol{v}}+\mathbf{h}\right).
\end{equation}
Here $\bm{h}$ satisfies the homogeneous boundary conditions of either
Dirichlet, Neumann or a combination of both, with boundary values
set to zero. The differential operator $\bm{\mathcal{L}}$ acting
on $\bm{h}$ now has to satisfy a homogeneous boundary value problem.
The overall strategy in using Fourier methods is to convert the differential
terms into spectral space using a spectral propagator matrix obtained
with a Fourier transform. The examples assume that the derivative
is second order, but the procedure is applicable to other more general
derivatives of the form given by Eq. (\ref{eq:linear-derivative}).

\subsection{Grid and transform methods}

In the examples we use an $M$-point time grid with a fixed time-step
$\Delta t$, and an $N$- point spatial grid with a fixed space step
$\Delta x$. Adaptive choices are also possible, and have advantages
in some cases \citep{mattar1980adaptive}. There is no unique choice
of trigonometric transform grid, although the spatial and momentum
grid spacings are restricted. Finite difference methods use a boundary
condition defined on a grid point, so we use discrete DST/DCT transforms
defined on whole grid points, as explained in Appendix C.

In all cases, if $N_{i}$ is the total number of grid points in the
$i-th$ dimension, the range is $R_{i}=X_{i}^{(b)}-X_{i}^{(a)}$ and
spatial step size $\Delta x_{i}=R_{i}/\left(N_{i}-1\right)$, so we
obtain a spatial grid for $j=1,\ldots N_{i}$ such that
\begin{equation}
\begin{split}x_{j}^{(i)}=X_{i}^{(a)}+(j-1)\Delta x_{i}.\end{split}
\end{equation}

The spatial indices used here start at $j=1$, but not all end points
are included in the usual definitions, since they are not needed if
the fields are zero. Discrete trigonometric transforms are closely
related to the fast Fourier transform or FFT \citep{cooley1965algorithm,Frigo1998,Frigo2005},
and have similar computational speed properties that scale proportional
to $N\ln N$, for $N$ spatial points in one dimension.

To compare with the more common FFT methods, if $N_{i}^{FT}$ is the
logical FFT size of a periodic discrete Fourier transform, then the
momenta are
\begin{equation}
\Delta k_{i}=\frac{2\pi}{N_{i}^{FT}\Delta x_{i}}.
\end{equation}
The corresponding momentum spacing in the case of a discrete sine
or cosine transform is:
\begin{equation}
\Delta k_{i}=\frac{\pi}{\left(N_{i}-1\right)\Delta x_{i}}=\frac{\pi}{R_{i}}.
\end{equation}

The momentum grid definitions depend on the transforms used, which
change with the boundary type $\beta$, as explained below. We use
standard terminology \citep{Frigo2005,britanak2010discrete} to identify
the different types of discrete transforms needed, depending on the
boundaries.

\subsubsection{Periodic grid:}

This has $N=N_{FT}$, and is included here for comparison purposes:
\begin{align}
k_{n} & =(n-1)\Delta k,\text{ for }~n=1,\ \ldots\ N
\end{align}

All spatial points with specified values are included in the FFT sum.

\subsubsection{Symmetric boundary: $\beta=1,4$}

In these cases, the boundary type is the same at the upper and lower
limits. Here, we use the DST-I and DCT-1 transforms, with:
\begin{align}
k_{n} & =(n-1)\Delta k,\text{ for }~n=1,\ \ldots\ N
\end{align}

\subsubsection{Non-symmetric boundary: $\beta=2,3$}

In these cases, the boundary type is different at the upper and lower
limits. In this case the DST-II/III and DCT-II/III transforms are
used, which have
\begin{align}
k_{n} & =(n-1/2)\Delta k,\text{ for }~n=1,\ \ldots\ N
\end{align}

\subsection{Fourier spectral derivatives}

In the Fourier spectral derivative method, the only term that is modified
is the derivative evaluated at a specific point in time. For clarity,
the derivative is evaluated here in the case of a Laplacian, but similar
methods hold in the more general case described in Eq. (\ref{eq:linear-derivative}).

Defining a discrete Fourier index $\bm{n}=$ $n_{1},...,n_{d}$, the
homogeneous variable $\mathbf{h}=\bm{u}-\bm{v}$ is transformed using
the appropriate DST or DCT, chosen according to the boundary types
$\beta_{i}^{\mu}$ for component $\mu$ and dimension $i$ as explained
in Appendix C, such that
\begin{align}
\tilde{h}_{\bm{n}}^{\mu} & =\sum_{\bm{j}=1}^{\bm{N}}h_{\bm{j}}^{\mu}\prod_{i=1}^{d}\mathcal{F}_{n_{i}j_{i}}^{\beta_{i}^{\mu}}\nonumber \\
h_{\bm{j}}^{\mu} & =\sum_{\bm{n}=1}^{\bm{N}}\tilde{h}_{\bm{n}}^{\mu}\prod_{i=1}^{d}\mathcal{I}_{j_{i}n_{i}}^{\beta_{i}^{\mu}},
\end{align}
where $\mathcal{F}$ is the forward and $\mathcal{I}$ the inverse
transform. The derivative term gives the result:
\begin{equation}
\bm{\mathcal{L}}\left[\boldsymbol{u}\right]^{\nu}=\sum_{\mu,o,i}D_{o,i}^{\nu,\mu}\sum_{\bm{n}=1}^{\bm{N}}k_{n_{i}}^{o}s_{o,i}^{\mu}\tilde{h}_{\bm{n}}^{\mu}\prod_{k=1}^{d}\mathcal{I}_{j_{k}n_{k}}^{\beta_{k}^{\mu}}+\bm{\mathcal{L}}\left[\boldsymbol{v}\right]^{\nu},\label{eq:Even-order case}
\end{equation}
where the inhomogeneous patch term $\boldsymbol{v}$ is typically
a low-order polynomial solution which is analytically differentiable.
The coefficient $s_{o,i}^{\mu}$ has unit magnitude, and depends on
the transform and the derivative order:

\subsubsection{Periodic case}

In periodic cases, if $\beta_{i}^{\mu}=0$, the coefficient is $s_{o,i}^{\mu}=i^{o}$,
where $i=\sqrt{-1}$, giving imaginary values for odd orders.

\subsubsection{Even-order case}

In non-periodic, even order cases, with $\beta_{i}^{\mu}>0$, the
coefficient is $s_{o,i}^{\mu}=\left(-1\right)^{o/2}$.

\subsubsection{Odd-order case}

In non-periodic, odd order cases, a sine changes to a cosine, and
vice-versa. The coefficient depends on the transform. Using the floor
function, $\left\lfloor .\right\rfloor $, one has that:
\begin{itemize}
\item $s_{o,i}^{\mu}=\left(-1\right)^{\left\lfloor o/2\right\rfloor }$
for DST, $\beta_{i}^{\mu}=1,4$
\item $s_{o,i}^{\mu}=\left(-1\right)^{\left\lfloor (o+1)/2\right\rfloor }$
for DCT, $\beta_{i}^{\mu}=2,3$
\end{itemize}

\subsubsection{Combined FSD expression}

Given the restrictions above, all the vector and grid indices can
be combined into a single index $\alpha$, to give a total of $C\prod_{i=1}^{d}N_{i}$
components. The field can be written as $V_{\alpha}$, and this transformation
leads to a discretized approximation to the derivatives,
\begin{equation}
\bm{\mathcal{L}}\left[\boldsymbol{v}\right]_{\alpha}\approx L_{\alpha\alpha'}V_{\alpha'}.
\end{equation}

Defining $H_{\alpha}$ as a single index form of $\mathbf{h}\left(\bm{r},\bm{v}\right)$,
the final PDE in the discretized space takes the form
\begin{equation}
\frac{\partial}{\partial t}V_{\alpha}=H_{\alpha}\left(\bm{r}_{\alpha},V_{\alpha}\right)+L_{\alpha\alpha'}V_{\alpha'}.
\end{equation}

This equation can be treated as an ordinary differential equation,
and solved in multiple ways. The examples treated here all use a midpoint
method in time. Higher order algorithms are possible. Both a forward
and backward transform are needed for each derivative evaluation.

\subsection{Fourier interaction picture}

An alternative strategy is a split-step or interaction picture method
\citep{BURSTEIN1970547,Hardin1973ApplicationOT,fisher1973role,Fisher1975,WERNER1997,caradoc2000vortex,Johan_JLT2007,Zhang:10,Balac2016}.
This approach has been very widely used in physics, especially for
treating propagation in nonlinear optics, quantum optics and Bose-Einstein
condensates. In this method, the linear part of the PDE, $\bm{\mathcal{L}}\cdot\boldsymbol{v}$
and the nonlinear part are separately integrated numerically, before
recombining the solutions. In previous methods, the linear part is
solved in the Fourier domain, with the assumption of periodic boundary
conditions.

We extend this to cases where $\bm{\mathcal{L}}$ is a linear derivative
term with constant coefficients, but with non-periodic boundaries,
while $\bm{h}\left[\bm{r},\bm{v}\right]$ includes linear or nonlinear
functions of the coordinates, derivatives and fields. We define the
interaction picture field $\bm{a}$ relative to a reference time $t'$
such that:
\begin{align}
v^{\mu}\left(t,\bm{x}\right) & =\int G^{\mu\nu}\left(\bm{x},t\left|\bm{x}',t'\right.\right)\cdot a^{\nu}\left(t,\bm{x}'\right)d\bm{x}'=\boldsymbol{\mathcal{\boldsymbol{G}}}\left[\bm{a}\right].
\end{align}
Here $\boldsymbol{G}$ is a Green's function defined such that:
\begin{itemize}
\item It obeys a linear partial differential equation:
\begin{equation}
\begin{split}\frac{\partial}{\partial t}\boldsymbol{G}\left(\bm{x},t\left|\bm{x}',t'\right.\right)=\bm{\mathcal{L}}\cdot\boldsymbol{G}\left(t,\bm{x}\left|\bm{x}'\right.\right).\end{split}
\end{equation}
\item It has an initial condition such that
\begin{equation}
G^{\mu\nu}\left(\bm{x},t'\left|\bm{x}',t'\right.\right)=\delta^{\mu\nu}\delta\left(\bm{x}-\bm{x}'\right).
\end{equation}
\item Boundary conditions are applied on the total field so that:
\begin{equation}
\mathcal{B}\left[\bm{r},\boldsymbol{\nabla},\boldsymbol{G}\left[\bm{a}\right]\right]=0.
\end{equation}
\end{itemize}
As a result, 
\begin{align}
\frac{\partial v^{\mu}}{\partial t}= & \int\left[\frac{\partial}{\partial t}G^{\mu\nu}\left(\bm{x},t\left|\bm{x}',t'\right.\right)\right]\cdot a^{\nu}\left(t,\bm{x}'\right)d\bm{x}'\nonumber \\
+ & \int G^{\mu\nu}\left(\bm{x},t\left|\bm{x}',t'\right.\right)\cdot\left[\frac{\partial}{\partial t}a^{\nu}\left(t,\bm{x}'\right)\right]d\bm{x}'.
\end{align}

In a more compact notation, omitting all time arguments for simplicity,
\begin{equation}
\frac{\partial\bm{v}}{\partial t}=\bm{\mathcal{L}}\cdot\bm{v}+\boldsymbol{G}\circ\left[\frac{\partial}{\partial t}\bm{a}\right].
\end{equation}
Therefore, provided that
\begin{align}
\boldsymbol{G}\circ\frac{\partial\bm{a}}{\partial t}= & \bm{h}\left[\bm{r},\bm{v}\right],
\end{align}
one has a solution to the original Eq. (\ref{eq:pde}). If the differential
operator has factorized solutions, $G$ is a diagonal matrix. These
are the examples treated here, but the method itself can be generalized.

This reduces the problem to an ordinary differential equation (ODE)
of form:
\begin{equation}
\frac{\partial\bm{a}}{\partial t}=\bm{\mathcal{A}}\left[\bm{a}\right],\label{eq:ODE}
\end{equation}
where $\bm{\mathcal{A}}$ is an array with $d+1$ indices, including
both internal and spatial degrees of freedom,

\subsection{Time-stepping algorithm}

Either the FSD or FIP method reduces to an ODE. There are many ways
to numerically solve an ODE of the form in Eq. (\ref{eq:ODE}), usually
with discretized time such that $t=t_{j}$, where $j=1,\ldots J$.
For clarity, we use a midpoint method (MP) in most of the examples,
although an RK4 time-stepping method is also possible \citep{caradoc2000vortex,Zhang:10},
and we give examples of this as well. 

The midpoint method is defined implicitly such that:
\begin{align}
\bm{a}\left(t_{j}+\Delta t/2\right) & =\bm{a}\left(t_{j}\right)+\Delta t\bm{\mathcal{A}}\left[\bm{a}\left(t_{j}+\Delta t/2\right)\right]/2\nonumber \\
\bm{a}\left(t_{j}+\Delta t\right) & =\bm{a}\left(t_{j}\right)+\Delta t\bm{\mathcal{A}}\left[\bm{a}\left(t_{j}+\Delta t/2\right)\right].
\end{align}
If we choose the interaction picture origin as $t'=\Delta t/2$, then
$\bm{a}=\bm{v}$ at the midpoint, and on redefining $\mathcal{\boldsymbol{G}}$
as a discretized version of the propagator, we obtain 
\begin{align}
\bm{a}\left(t_{j}\right) & =\boldsymbol{\mathcal{\boldsymbol{G}}}^{-1}\left(-\Delta t/2\right)\left[\bm{v}\left(t_{j}\right)\right]\nonumber \\
\bm{v}\left(t_{j+1}\right) & =\boldsymbol{\mathcal{\boldsymbol{G}}}\left(\Delta t/2\right)\left[\bm{a}\left(t_{j+1}\right)\right].
\end{align}
 In terms of the original field, this result becomes: 
\begin{equation}
\bm{v}\left(t_{j+1}\right)=\boldsymbol{\mathcal{\boldsymbol{G}}}\left(\Delta t\right)\left[\bm{v}\left(t_{j}\right)\right]+\Delta t\boldsymbol{\mathcal{\boldsymbol{G}}}\left(\Delta t/2\right)\mathcal{\bm{\mathcal{A}}}\left[\bm{v}\left(t_{j}+\Delta/2\right)\right].
\end{equation}
At the midpoint, the derivative $\mathcal{\bm{\mathcal{A}}}\left[\bm{v}\right]$
is found iteratively, since this is an implicit method. Similar techniques
can be used for Runge-Kutta methods.

\subsection{Non-periodic boundaries}

To define the numerical method one must also specify the boundary
values for $\boldsymbol{\mathcal{\boldsymbol{G}}}\left(t\right)$,
which may only be available at discrete lattice points. We impose
boundary values that solve the boundary-value equation at the lattice
points and are constant for $t_{j}<t<t_{j}+\Delta t/2$ and for $t_{j}+\Delta t/2<t<t_{j+1}$
. Hence, $\boldsymbol{\mathcal{\boldsymbol{G}}}^{-1}\left(-\Delta t/2\right)$
is defined with respect to the initial boundary values, and $\boldsymbol{\mathcal{\boldsymbol{G}}}\left(\Delta t/2\right)$
is defined with respect to the final boundary values, so that:
\begin{equation}
\mathcal{B}\left[\bm{r},\boldsymbol{\nabla},\bm{a}\left(t_{k}\right)\right]=0.
\end{equation}

We treat the Green's function part of the problem using a modified
tau method, in which the intermediate solution 
\begin{align}
\bm{u}_{I} & =\boldsymbol{\mathcal{\boldsymbol{G}}}\left[\bm{a}\right]\nonumber \\
 & =\bm{h}+\bm{v}\label{eq:two parts}
\end{align}
 is divided up into a homogeneous ($h$) and inhomogeneous ($v$)
part. In the interaction picture for multiple derivatives, the derivatives
are assumed to be a sum of independent terms for each dimension and
field component, which are solved in parallel, such that:
\begin{equation}
\bm{\mathcal{L}}\left[\boldsymbol{\nabla}\right]\cdot\boldsymbol{v}=\sum_{i}\left[\mathcal{L}_{x}+\mathcal{L}_{y}+\mathcal{L}_{z}+...\right]\bm{v}=\dot{\bm{v}}
\end{equation}
and the inhomogeneous part satisfies the boundary conditions given
by 
\begin{equation}
\mathcal{B}\left[\bm{r},\boldsymbol{\nabla},\bm{v}\left(t_{k}\right)\right]=0.
\end{equation}

This requires the existence of initial conditions for which the inhomogeneous
part is readily solvable, which depends on the form of the boundary
values. For clarity in the description, we assume that there are only
second derivative terms, and hence one can treat each dimension separately:
\begin{equation}
\mathcal{L}_{i}a_{i}=D_{i}^{(2)}\frac{\partial^{2}}{\partial x^{2}}a_{i}.
\end{equation}
We initially treat only one space dimension, but this approach can
be generalized to arbitrary numbers of dimensions, again depending
on the boundary values. On dividing up the solution into homogenous
and inhomogeneous parts as in Eq.(\ref{eq:two parts}), for a single
dimension, each component satisfies the Laplacian equation
\begin{equation}
\frac{\partial\left(\bm{h},\bm{v}\right)}{\partial t}=D^{(2)}\boldsymbol{\nabla}^{2}\left(\bm{h},\bm{v}\right).
\end{equation}

The boundary terms are applied with zero boundaries for the homogenous
part. There can be initial conditions applied to each of $h$ and
$v$. The initial conditions applied to $h$ must be either periodic,
or else give zero boundary values for either the field $h$ (Dirichlet
case) or its derivative (Robin/Neumann case) at the upper and lower
spatial boundaries. In the case of the inhomogeneous part, the boundary
values of either the $v$ field or its derivative are considered to
be nonzero.

\section{Inhomogeneous patch solutions\label{sec:Inhomogeneous-patch-solutions}}

For inhomogeneous boundary value problems, this approach involves
adding a patch function \citep{boyd2001}, typically a polynomial
function. This patch function is smooth and must satisfy the boundary
conditions. In the FIP case it must locally solve the differential
equation over a short time. Here we explain how a patch function can
be selected to meet the boundary conditions of the problem.

Just as there are different homogeneous cases, there are four cases
for the inhomogeneous part, which require a combination of solutions
of two of $\bm{U}_{a,b}$ (Dirichlet) and $\bm{N}_{a,b}$ (Neumann).
The prefactors are given explicitly as $\bm{U}_{a,b}\left(\bm{u}\left(t_{k}\right),t_{k}\right)$
and $\bm{N}_{a,b}\left(\bm{u}\left(t_{k}\right),t_{k}\right)$. The
post factors, in general, are known only implicitly as $\bm{U}_{a,b}\left(\bm{u}\left(t_{k+1}\right),t_{k+1}\right)$
and $\bm{N}_{a,b}\left(\bm{u}\left(t_{k+1}\right),t_{k+1}\right)$.
This may require them to be treated iteratively if the boundary value
depends on $\bm{u}$ as well as $t$.

We assume a single, diagonal second-order derivative. Any higher order
derivative terms all vanish exactly for these patches. Also, for simplicity
we suppose that the boundary condition is approximately time-independent
over a short time-interval.

\subsection{Boundary types}

The inhomogeneous patch function $\bm{v}$ is obtained for the four
combinations of boundary types, in one space dimension. This can be
extended to higher dimensions if the boundary functions factorize,
but otherwise more complex patches are needed.

\subsubsection{Dirichlet-Dirichlet}

The boundary condition of the problem is defined such that
\begin{align}
\bm{v}\left(0,x_{a,b}\right) & =\bm{U}_{a,b}.
\end{align}
We define a patch function of the form 
\begin{equation}
\bm{v}\left(t,x\right)=\bm{U}_{a}+\frac{x-x_{a}}{x_{b}-x_{a}}\left(\bm{U}_{b}-\bm{U}_{a}\right)
\end{equation}
which satisfies the condition
\begin{equation}
\begin{split}\frac{\partial}{\partial t}\left[\bm{v}\right]=D^{(2)}\boldsymbol{\nabla}^{2}\cdot\left[\bm{v}\right]=0.\end{split}
\end{equation}

\subsubsection{Dirichlet-Neumann:}

When considering a combination of boundary conditions such that

\begin{align}
\bm{v}(0,x_{a}) & =\bm{U}_{a}\nonumber \\
\partial_{x}\bm{v}(0,x_{a}) & =\bm{N}_{b},
\end{align}
we define a patch function of the form
\begin{equation}
\bm{v}\left(t,x\right)=\bm{U}_{a}+\left(x-x_{a}\right)n_{b},
\end{equation}
which satisfies
\begin{equation}
\begin{split}\frac{\partial}{\partial t}\left[\bm{v}\right]=D^{(2)}\boldsymbol{\nabla}^{2}\cdot\left[\bm{v}\right]=0.\end{split}
\end{equation}

\subsubsection{Neumann-Dirichlet:}

Similarly, if the boundary conditions are defined such that

\begin{align}
\bm{v}(0,x_{b}) & =\bm{U}_{b}\nonumber \\
\partial_{x}\bm{v}(0,x_{b}) & =\bm{N}_{a},
\end{align}
we define a patch function of the form
\begin{equation}
\bm{v}\left(t,x\right)=\bm{U}_{b}+\left(x-x_{b}\right)\bm{N}_{a},
\end{equation}
which satisfies
\begin{equation}
\begin{split}\frac{\partial}{\partial t}\left[\bm{v}\right]=D_{2}\boldsymbol{\nabla}^{2}\cdot\left[\bm{v}\right]=0.\end{split}
\end{equation}

\subsubsection{Neumann-Neumann:}

The boundary condition of the problem is defined such that:

\begin{equation}
\partial_{x}\bm{v}\left(t_{0},x_{a,b}\right)=\bm{N}_{a,b}.
\end{equation}
We define a patch function with space-derivative
\begin{equation}
\partial_{x}\bm{v}\left(t_{0},x\right)=\bm{N}_{a}+\frac{x-x_{a}}{x_{b}-x_{a}}\left(\bm{N}_{b}-\bm{N}_{a}\right),
\end{equation}
which leads to
\begin{align}
\bm{v}\left(t,x\right) & =\bm{\epsilon}\left(t-t_{0}\right)+\bm{N}_{a}\left(x-x_{a}\right)\nonumber \\
 & \ +\frac{1}{2}\frac{\left(x-x_{a}\right)^{2}}{x_{b}-x_{a}}\left(\bm{N}_{b}-\bm{N}_{a}\right).
\end{align}

This satisfies
\begin{align}
\frac{\partial}{\partial t}\left[\bm{v}\right] & =\bm{\epsilon}\\
D_{2}\nabla^{2}\left[\bm{v}\right] & =\frac{D_{2}\left(\bm{N}_{b}-\bm{N}_{a}\right)}{\left(x_{b}-x_{a}\right)}
\end{align}
Hence, 
\begin{equation}
\bm{\epsilon}=\frac{D_{2}\left(\bm{N}_{b}-\bm{N}_{a}\right)\Delta t}{\left(x_{b}-x_{a}\right)}
\end{equation}

\subsection{Algorithm \label{subsec:Algorithm}}

Both the Fourier spectral derivative (FSD) and the Fourier interaction
picture (FIP) methods numerically solve the PDE in Eq. (\ref{eq:pde})
by using the midpoint algorithm to propagate the field $u$ in time
in most examples. The method is explained here in detail, and can
be extended to a fourth-order Runge-Kutta algorithm for higher-order
convergence in time.

The differences between FSD and FIP when calling the midpoint algorithm
lie in how the propagation function $\mathcal{D}$ and the linear
propagator $\mathcal{P}$ are defined. In FSD, we apply the transforms
detailed in previous sections to the propagation function $\mathcal{D}\equiv g\left[\bm{r},\boldsymbol{\nabla},\bm{u}\right]+\bm{\mathcal{L}}\left[\boldsymbol{\nabla}\right]\cdot\boldsymbol{u}$
and the linear propagator is just $1$. In FIP, these transforms are
only applied to the $g$ term, hence defining $\mathcal{D}\equiv g\left[\bm{r},\boldsymbol{\nabla},\bm{u}\right]$,
while the linear operator $\bm{\mathcal{L}}\left[\boldsymbol{\nabla}\right]$
term is used to define the linear propagator $\mathcal{P}\left(t',t\right)=\text{exp}\left[\left(t'-t\right)\mathcal{L}[\nabla]\right]$.

The steps involved in the midpoint algorithm are outlined in the pseudo-code
below:\\

\noindent\doublebox{\begin{minipage}[t]{1\columnwidth - 2\fboxsep - 7.5\fboxrule - 1pt}%
\textbf{MIDPOINT ALGORITHM}

\medskip{}

$\bar{u}^{(0)}=\mathcal{P}\left(t_{j}+\frac{\Delta t}{2},t_{j}\right)u_{j}$

\smallskip{}

\textbf{for }$i$ from $1$ to $iter$ \textbf{do}

\smallskip{}

~~~~~~~~~~~~$\bar{u}^{(i)}=\bar{u}^{(0)}+\frac{\Delta t}{2}\,\mathcal{D}\left[\bar{u}^{(i-1)},t_{j}+\frac{\Delta t}{2}\right]$

\smallskip{}
\textbf{end}

\smallskip{}

\textbf{$u_{j+1}=\mathcal{P}\left(t_{j}+\Delta t,t_{j}+\frac{\Delta t}{2}\right)\left[2\bar{u}^{(iter)}-\bar{u}^{(0)}\right]$}%
\end{minipage}}

\section{Comparisons between uniform and non-uniform grid spacing\label{sec:Comparisons-between-uniform}}

Here, we compare the errors of the FSD and FIP Fourier spectral methods
with those of other PDE solvers that use polynomial spectral methods.
The FIP and FSD methods in all examples were implemented in xSPDE4
\citep{Kiesewetter2016,Kiesewetter2023,Drummond2023}, a public domain
code in the Matlab language. Fourier transforms were implemented with
the Matlab interface to FFTW functions, a public domain fast Fourier
transform package written in C. The DCT and DST transforms are explained
in Appendix C, and were obtained by computing complex Fourier transforms
after using symmetries. Further efficiency improvements are possible
with dedicated DCT/DST code.

We consider two public domain PDE solvers: (1) Dedalus, which is a
public-domain package that uses spectral polynomial methods with unequal
spaced grids \citep{Burns2020} and (2) Pdepe, which has a public
domain implementation that is now an inbuilt function in MATLAB language,
using a method of lines approach with an equally spaced grid \citep{Skeel1990}.
We consider PDE examples with analytical solutions that allow errors
to be calculated. These solutions consist of both relatively smooth
and rapidly varying functions of space and time in order to properly
compare and contrast the algorithms. The comparisons are focused on
determining accuracy and speed with fixed numbers of space-time steps.
The equation for the relevant RMS relative comparison error used for
comparisons is discussed in Appendix D.

We consider up to four combinations of boundary conditions for each
field component and dimension. These are the Dirichlet-Dirichlet (D-D),
Dirichlet-Neumann(D-N), Neumann-Dirichlet (N-D), and Neumann-Neumann
(N-N). Most boundary conditions used are time-dependent to better
align with cases of practical interest. The errors reported are the
RMS relative comparison error, normalized by the maximum of each computed
output. The Fourier results obtained use either the FIP or the FSD
algorithm with the midpoint method (MP) as the integration method.
Since the main focus is to evaluate the influence of using discrete
trigonometric transforms on spatial discretization error, we use a
simple time-stepping method (midpoint algorithm) for the following
comparisons. 

The Tau results were obtained using the Dedalus software, employing
either a Chebyshev or Legendre basis, an RK222 integration method
\citep{Ascher1997implicit} and a dealiasing factor of $d=1$ or $d=2$.
 All problems were solved using both Chebyshev and Legendre bases,
and in each case we report the results corresponding to the basis
that yields the lowest error. The Galerkin results used the Matlab
pdepe software, with an absolute and relative error tolerance set
to 1, and the InitialStep and MaxStep options for the ode solver set
equal to the time step-size $\Delta t$. 

All results are computed using a 2020 MacBook Pro, to eliminate hardware
differences, with Fourier results highlighted in red.

\subsection{Linear PDE example}

\subsubsection{Heat equation with non-periodic boundaries}

This example solves a (1+1)-dimensional PDE with an initial condition
of $\mathbf{u}\left(t=0,x\right)=\mathbf{f}\left(x\right)$ and

\begin{eqnarray}
\frac{\partial\mathbf{u}}{\partial t} & = & \frac{\partial^{2}\mathbf{u}}{\partial x^{2}}\,.
\end{eqnarray}
The solution is subject to either Dirichlet and/or Neumann boundary
conditions with boundary values of zero at $x_{\pm}=[0,\pi]$ so that
$u\left(t,x_{\pm}\right)=0\,$ or $\partial u/\partial x\left(t,x_{\pm}\right)=0$.
Each component has different combinations of boundary types. In all
cases the grid range is from $x=0$ to $x=\pi$, and the time duration
is from $t=0$ to $t=4$.

\paragraph*{1. Dirichlet-Dirichlet}

With $u(0)=u(\pi)=0$, the exact solution has the form: 
\begin{align}
u & =\sum_{n=1}^{2}S_{n}\sin\left(nx\right)e^{-n^{2}t}.
\end{align}

For this case, a solution is: 
\begin{equation}
u(x,t)=4\sin\left(x\right)e^{-t}+\sin\left(2x\right)e^{-4t}.
\end{equation}

\paragraph*{2. Neumann-Neumann}

With $\partial_{x}u(0)=\partial_{x}u(\pi)=0$, the exact solution
has the form: 
\begin{align}
u & =\sum_{n=0}^{\infty}C_{n}\cos\left(nx\right)e^{-n^{2}t}.
\end{align}

For this case, a solution is: 
\begin{equation}
u(x,t)=5+4\cos\left(x\right)e^{-t}+\cos\left(2x\right)e^{-4t}.
\end{equation}

\paragraph*{3. Dirichlet-Neumann}

Here $u(0)=\partial_{x}u(\pi)=0$, the exact solution has the form:
\begin{align}
u & =\sum_{n=1}^{\infty}S_{n}\sin\left((2n-1)x/2\right)e^{-(2n-1)^{2}t/4}.
\end{align}

For this case, a solution is: 
\begin{equation}
u(x,0)=4\sin\left(x/2\right)e^{-t/4}+\sin\left(3x/2\right)e^{-9t/4}.
\end{equation}

\paragraph*{4. Neumann-Dirichlet}

Here $\partial_{x}u(0)=u(\pi)=0$, the general solution has the form:
\begin{align}
u & =\sum_{n=1}^{\infty}C_{n}\cos\left((2n-1)x/2\right)e^{-(2n-1)^{2}t/4}.
\end{align}

For this case, a solution is: 
\begin{equation}
u(x,t)=4\cos\left(x/2\right)e^{-t/4}+\cos\left(3x/2\right)e^{-9t/4}.
\end{equation}

\begin{figure}[H]
\centering{}\includegraphics[width=0.9\columnwidth]{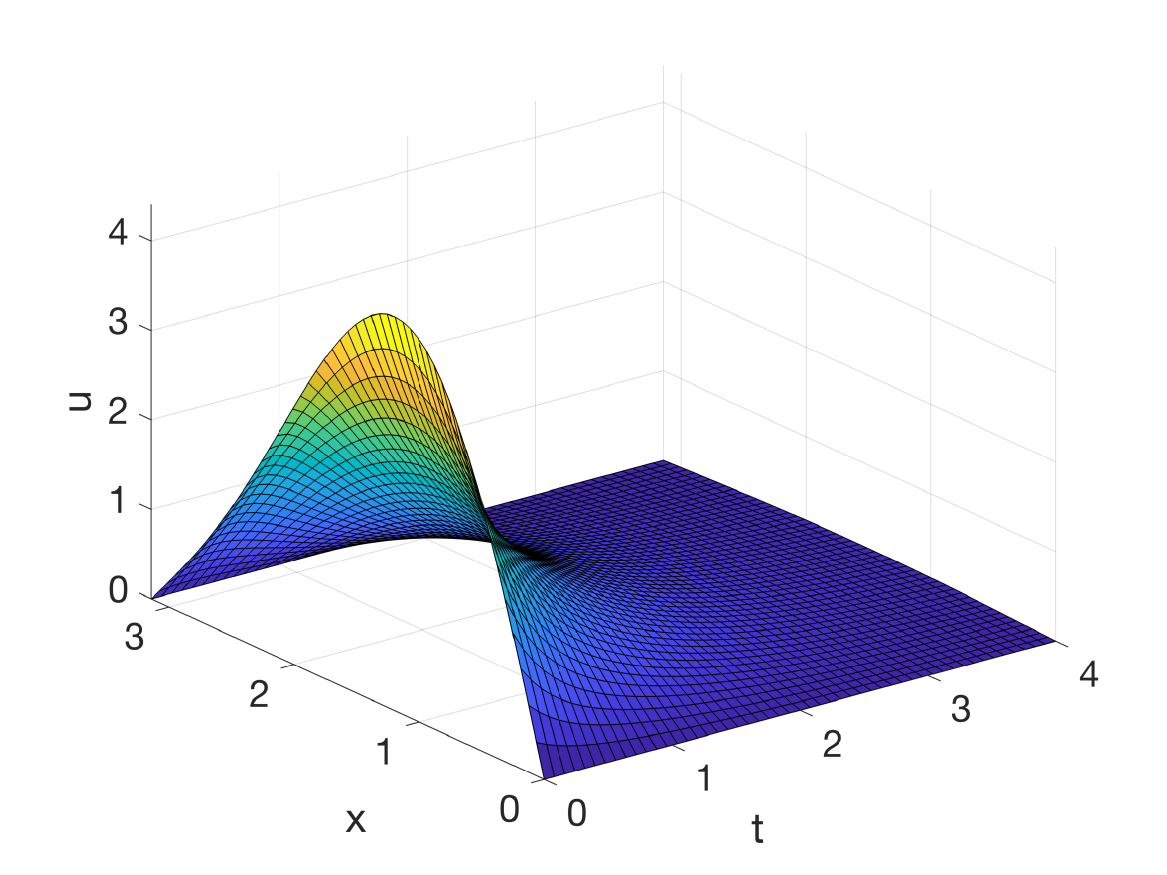}\caption{Plot of the solution to the heat equation using the Fourier interaction
picture (FIP) method for the D-D boundary condition with $\Delta x=\pi/50$
and $\Delta t=0.08$. \label{fig:heatFIP}}
\end{figure}
\begin{table}[!h]
\begin{centering}
\begin{tabular}{|c|c|c|c|c|c|c|}
\hline 
Bdary & FIP error & Tc (s) & Tau error & Tc (s) & Gal error & Tc (s)\tabularnewline
\hline 
\hline 
D-D & \textcolor{red}{$\mathbf{3\times10^{-16}}$} & \textbf{\textcolor{red}{$\mathbf{0.03}$}} & $4.2\times10^{-4}$ & $0.08$ & $1.1\times10^{-4}$ & $0.27$\tabularnewline
\hline 
N-N & \textcolor{red}{$\mathbf{2\times10^{-16}}$} & \textbf{\textcolor{red}{$\mathbf{0.03}$}} & $4.2\times10^{-4}$ & $0.08$ & $4.8\times10^{-5}$ & $0.17$\tabularnewline
\hline 
D-N & \textcolor{red}{$\mathbf{6\times10^{-15}}$} & \textbf{\textcolor{red}{$\mathbf{0.03}$}} & $1.6\times10^{-4}$ & $0.08$ & $4.6\times10^{-5}$ & $0.22$\tabularnewline
\hline 
N-D & \textcolor{red}{$\mathbf{2\times10^{-15}}$} & \textbf{\textcolor{red}{$\mathbf{0.03}$}} & $1.6\times10^{-4}$ & $0.08$ & $3.3\times10^{-5}$ & $0.29$\tabularnewline
\hline 
\end{tabular}
\par\end{centering}
\centering{}\caption{Comparison of the errors obtained for solving the heat equation using
the \textcolor{black}{Fourier inter}action picture (FIP), Tau and
Galerkin (Gal) methods with different boundary conditions. The number
of spatial and time steps is 50, with $\Delta x=\pi/50$ and $\Delta t=0.08$.
The Tau results were obtained using a Chebyshev basis with dealiasing
factor $d=1$.\label{tab:Heat_TD}}
\end{table}

In this example of the linear heat equation, the errors of the FIP
method are close to the machine precision limit (Table \ref{tab:Heat_TD}).
The Galerkin method produced smaller errors than the Tau method. However,
both of these have more than $10^{12}$ times greater errors than
the FIP algorithm. The run-time of the Galerkin method was approximately
10 times greater, and that of the Tau method about 3 times greater
than the FIP method. For this type of linear problem, the Fourier
interaction picture method appears to be almost ideally suited.

\subsection{Nonlinear PDE examples}

\subsubsection{Peregrine solitary wave with arbitrary boundary conditions \label{subsec:Peregrine-solitary-wave}}

Peregrine solitary waves are models for isolated large ocean waves
\citep{Peregrine1983}. They are modeled as solutions to a (1+1)-dimensional
PDE, the nonlinear Schrödinger equation
\begin{eqnarray}
\frac{\partial u}{\partial t} & = & i\cdot\left(u\cdot\left|u\right|^{2}+\frac{1}{2}\frac{\partial^{2}u}{\partial x^{2}}\right)\,.
\end{eqnarray}
The Peregrine solution on an infinite domain is:
\begin{eqnarray}
Per(t,x) & = & e^{it}\left(\frac{4\left(1+2it\right)}{1+4\left(t^{2}+x^{2}\right)}-1\right)\,.
\end{eqnarray}

In the example, this is solved using finite boundary conditions, with
initial and boundary values corresponding to the exact solution. The
solution is subject to four combinations of Dirichlet or Neumann boundary
conditions with boundary values at $x_{m}=\pm2$, and four different
combinations of $u\left(t,\pm2\right)=Per(t,\pm2)$, and $\frac{\partial u}{\partial x}\left(t,\pm2\right)=\left.\frac{\partial}{\partial x}Per(t,x)\right|_{\pm2}\,.$
The grid range is from $x=-2$ to $x=2$, and the time duration is
from $t=-5$ to $t=5$. The observable is $o\equiv\left|u\right|^{2}$.

\begin{figure}[H]
\centering{}\includegraphics[width=0.9\columnwidth]{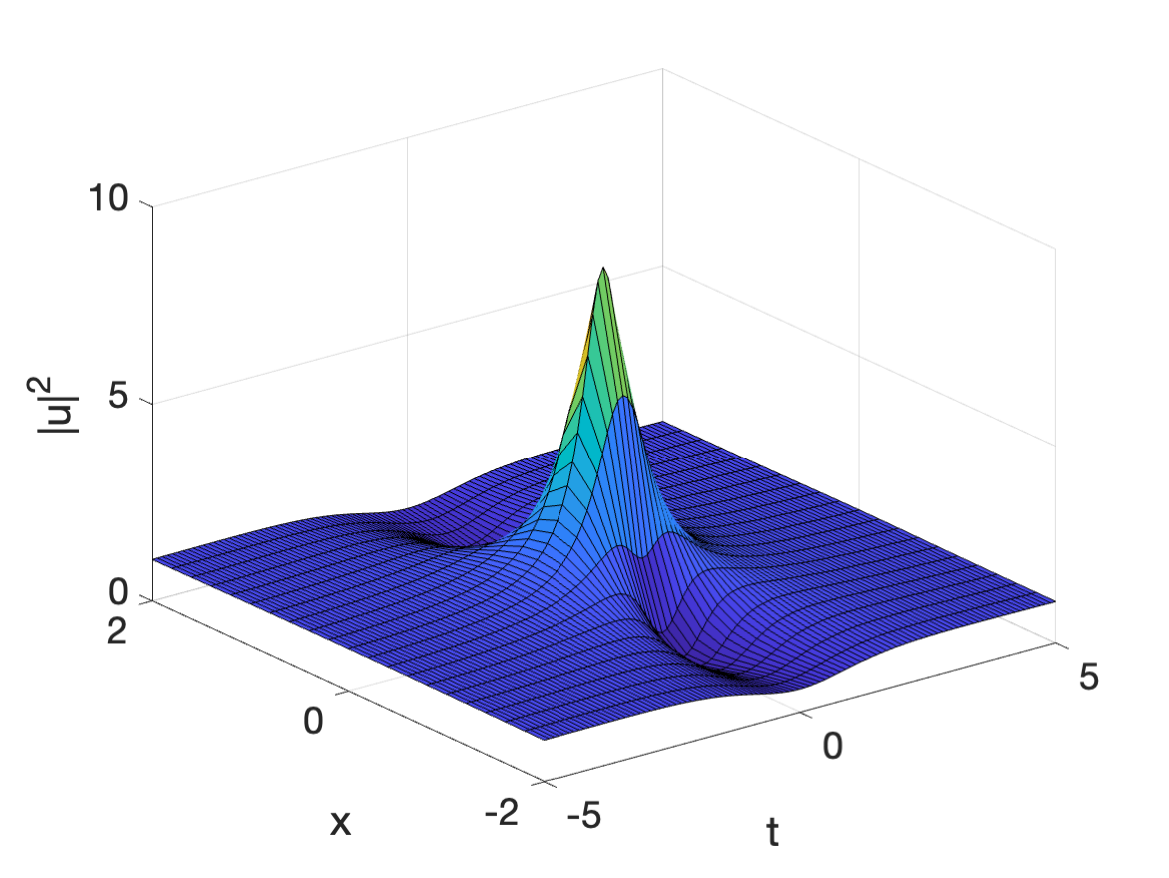}\caption{Plots of the solution to the Peregrine solitary wave equation using
the Fourier spectral derivative (FSD) method with arbitrary boundaries
for D-D boundary condition with $\Delta x=0.2$ and $\Delta t=0.005$.
The other graphs are nearly identical.\label{fig:Peregrine}}
\end{figure}

\begin{table}[!h]
\begin{centering}
\begin{tabular}{|c|c|c|c|c|}
\hline 
Boundary & FSD error & Tc (s) & FIP error & Tc( s)\tabularnewline
\hline 
\hline 
D-D & \textcolor{red}{$\mathbf{3.3\times10^{-4}}$} & \textcolor{red}{$\mathbf{0.9}$} & \textcolor{red}{$\mathbf{3.3\times10^{-4}}$} & \textcolor{red}{$\mathbf{0.8}$}\tabularnewline
\hline 
N-N & \textcolor{red}{$\mathbf{3\times10^{-4}}$} & \textcolor{red}{$\mathbf{1.0}$} & \textcolor{red}{$\mathbf{1.1\times10^{-3}}$} & \textcolor{red}{$\mathbf{0.7}$}\tabularnewline
\hline 
D-N & \textcolor{red}{$\mathbf{1\times10^{-3}}$} & \textcolor{red}{$\mathbf{1.4}$} & \textcolor{red}{$\mathbf{2.6\times10^{-3}}$} & \textcolor{red}{$\mathbf{0.9}$}\tabularnewline
\hline 
N-D & \textcolor{red}{$\mathbf{1\times10^{-3}}$} & \textcolor{red}{$\mathbf{1.3}$} & \textcolor{red}{$\mathbf{2.6\times10^{-3}}$} & \textcolor{red}{$\mathbf{0.9}$}\tabularnewline
\hline 
\end{tabular}\medskip{}
\par\end{centering}
\begin{centering}
\begin{tabular}{|c|c|c|c|c|}
\hline 
Boundary & Tau error & Tc (s) & Gal error & Tc (s)\tabularnewline
\hline 
\hline 
D-D & $0.01$ & $2.8$ & $0.23$ & $7.1$\tabularnewline
\hline 
N-N & $0.02$ & $2.8$ & $0.67$ & $3.7$\tabularnewline
\hline 
D-N & $0.03$ & $2.8$ & $0.52$ & $8.4$\tabularnewline
\hline 
N-D & $0.03$ & $2.8$ & $0.52$ & $2.1$\tabularnewline
\hline 
\end{tabular}
\par\end{centering}
\centering{}\caption{Comparison of errors obtained from solving the Peregrine solitary
wave equation using the Fourier spectral derivative (FSD), Fourier
interaction picture (FIP), Tau and Galerkin (Gal) method with different
boundary conditions. The number of spatial steps is 20 and time steps
is 2000 with $\Delta x=0.2$ and$\Delta t=0.005$. The Tau results
were obtained using a Chebyshev basis with dealiasing factor $d=2$.
All the errors are normalized by the maximum value $m=9.$\label{tab:Peregrine}}
\end{table}

The Peregrine solution is an isolated peak in space-time, causing
rapid variation in space near the peak. Based on the results of Table
\ref{tab:Peregrine}, this is handled reasonably accurately with FSD
and FIP methods in all four types of boundary condition treated here.
By comparison, the Tau and Gal methods experience errors for these
parameter values of nearly two and three orders of magnitude more,
respectively. The Galerkin method gave especially large errors, and
generally, such rapid nonlinear variations are most efficiently treated
using the Fourier spectral method.

\subsubsection{Breather}

Breathers of the nonlinear Schrödinger equation (NLSE) are specific
types of localized solutions with periodic behavior in both space
and time. They are solutions to a (1+1)-dimensional PDE,

\begin{equation}
\frac{\partial u}{\partial t}=-i\left(u|u|^{2}+\frac{1}{2}\frac{\partial^{2}u}{\partial x^{2}}\right).
\end{equation}

We consider a particular solution as mentioned in \citep{Satsuma1974}
with an initial condition of $u\left(t=0,x\right)=2sech\left(x\right)$,
where the envelope pulsates at a frequency of $\pi/2$. The solution
is subject to four combinations of Dirichlet or Neumann boundary conditions
with boundary values of four different combinations at the boundary
$\pm x_{m}$ with

\begin{equation}
u(t,x)=4e^{-it/2}\left(\frac{cosh(3x)+3e^{-4it}cosh(x)}{cosh(4x)+4cosh(2x)+3cos(4t)}\right)
\end{equation}

The solutions of this type are examined in optical fibers \citep{Gordon1983},
plasma physics and Bose-Einstein condensates \citep{Luo2020}. The
grid range is from $x=-2$ to $x=2$, and the time duration is from
$t=0$ to $t=\pi$. The observable is $o\equiv|u|$.

\begin{figure}[H]
\centering{}\includegraphics[width=0.9\columnwidth]{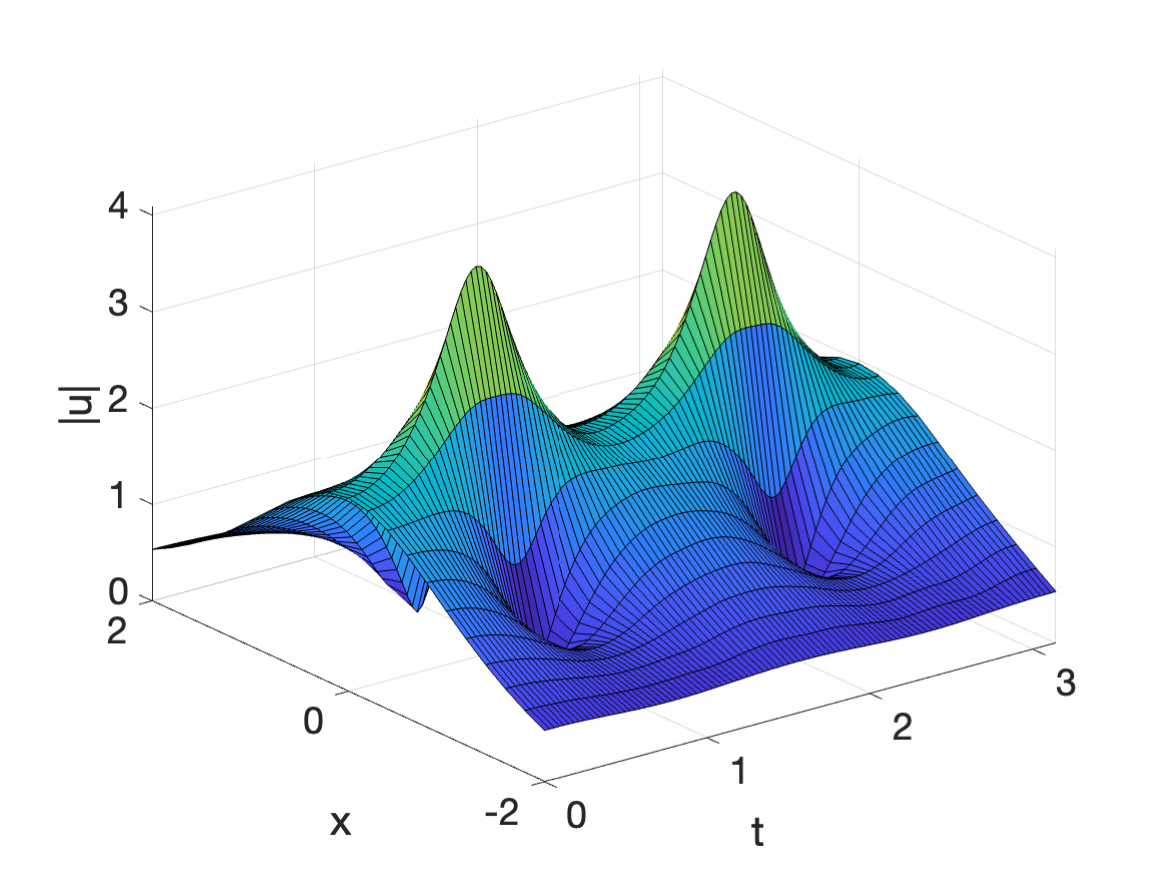}\caption{Plots of the breather solution of the nonlinear Schrödinger equation
using the Fourier spectral derivative (FSD) method with arbitrary
boundaries for the D-D boundary condition with $\Delta x=0.2$ and
$\Delta t=\pi/2000$.\label{fig:Breather}}
\end{figure}

\begin{table}[!h]
\begin{centering}
\begin{tabular}{|c|c|c|c|c|}
\hline 
Boundary & FSD error & Tc (s) & FIP error & Tc (s)\tabularnewline
\hline 
\hline 
D-D & \textcolor{red}{$\mathbf{5.03\times10^{-3}}$} & \textcolor{red}{$\mathbf{0.98}$} & \textcolor{red}{$\mathbf{4.95\times10^{-3}}$} & \textcolor{red}{$\mathbf{0.71}$}\tabularnewline
\hline 
N-N & \textcolor{red}{$\mathbf{4.38\times10^{-3}}$} & \textcolor{red}{$\mathbf{0.92}$} & \textcolor{red}{$\mathbf{4.31\times10^{-3}}$} & \textcolor{red}{$\mathbf{0.65}$}\tabularnewline
\hline 
D-N & \textcolor{red}{$\mathbf{5.63\times10^{-3}}$} & \textcolor{red}{$\mathbf{1.39}$} & \textcolor{red}{$\mathbf{5.54\times10^{-3}}$} & \textcolor{red}{$\mathbf{0.91}$}\tabularnewline
\hline 
N-D & \textcolor{red}{$\mathbf{5.63\times10^{-3}}$} & \textcolor{red}{$\mathbf{1.17}$} & \textcolor{red}{$\mathbf{5.54\times10^{-3}}$} & \textcolor{red}{$\mathbf{0.85}$}\tabularnewline
\hline 
\end{tabular}
\par\end{centering}
\begin{centering}
\medskip{}
\begin{tabular}{|c|c|c|c|c|}
\hline 
Boundary & Tau error & Tc (s) & Gal error & Tc (s)\tabularnewline
\hline 
\hline 
D-D & $0.05$ & $2.4$ & $0.11$ & $7.8$\tabularnewline
\hline 
N-N & $0.05$ & $2.3$ & $0.13$ & $7.5$\tabularnewline
\hline 
D-N & $0.08$ & $2.3$ & $0.10$ & $7.8$\tabularnewline
\hline 
N-D & $0.08$ & $2.3$ & $0.14$ & $7.6$\tabularnewline
\hline 
\end{tabular}
\par\end{centering}
\centering{}\caption{Comparison of errors obtained from solving the breather solution of
the nonlinear Schrödinger equation using the Fourier spectral derivative
(FSD), Fourier interaction picture (FIP), Tau and Galerkin (Gal) method
with different boundary conditions. The number of spatial steps is
20, and time steps is 2000 with $\Delta x=0.2$ and $\Delta t=\pi/2000$.
The Tau results were obtained using a Chebyshev basis with dealiasing
factor $d=2$.\label{tab:Breather}}
\end{table}

From Table \ref{tab:Breather}, we see that the very rapid spatial
variation occurring near the breather focus has relatively little
effect on the FSD and FIP error, but causes an error increase of nearly
one order of magnitude in the Tau and two orders of magnitude for
the Galerkin methods used for comparison purposes.

\subsubsection{Double simultons with phase variation}

Simultons occur in parametric waveguides \citep{Werner1993,Werner1997_simulton,He1998},
with the equation:
\begin{eqnarray}
\frac{\partial u_{1}}{\partial t} & = & -i\left[\frac{\partial^{2}u_{1}}{\partial x^{2}}+u_{1}^{*}u_{2}\right]\\
\frac{\partial u_{2}}{\partial t} & = & -i\left[\frac{\partial^{2}u_{2}}{\partial x^{2}}+u_{1}^{2}+u_{2}\right].\nonumber 
\end{eqnarray}
The simulton in one space dimension is

\begin{equation}
[u_{1},u_{2}](t,x)=\frac{3}{2}sech^{2}\left(\frac{x}{2}\right)[e^{-it},e^{-2it}]
\end{equation}

with derivatives
\begin{equation}
\frac{\partial[u_{1},u_{2}](t,\pm x)}{\partial x}=\mp\frac{3}{2}sech^{2}\left(\frac{x}{2}\right)tanh\left(\frac{x}{2}\right)[e^{-it},e^{-2it}]\,.\,\;
\end{equation}

The solution is subject to four combinations of Dirichlet and Neumann
boundary conditions at $x_{m}=\pm3$, where the boundaries are specified
differently for each component, and where:
\begin{eqnarray}
[u_{1},u_{2}]\left(t,\pm x_{m}\right) & = & [U_{1},U_{2}](t,\pm x_{m})\nonumber \\
\frac{\partial[u_{1},u_{2}]}{\partial x}\left(t,\pm x_{m}\right) & = & \left.\frac{\partial}{\partial x}[U_{1},U_{2}](t,x)\right|_{\pm x_{m}}\,.
\end{eqnarray}

The grid range is from $x=-3$ to $x=3$, and the time duration is
from $t=0$ to $t=\pi$. The observables are $o_{1}\equiv\Re\left(u_{1}\right)$
, and $o_{2}\equiv\Re\left(u_{2}\right)$ , giving two sets of four
outputs.

\begin{figure}[H]
\begin{centering}
\includegraphics[width=0.9\columnwidth,height=5cm]{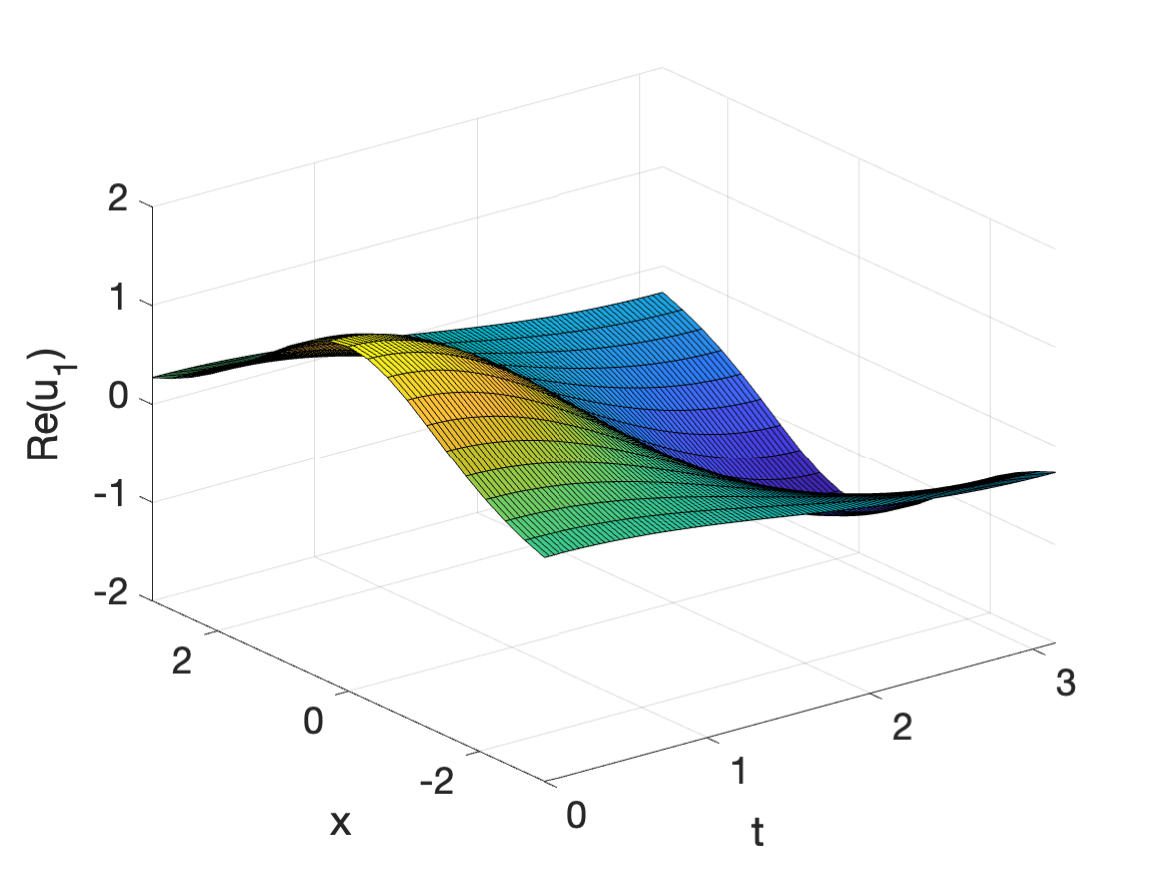}
\par\end{centering}
\begin{centering}
\includegraphics[width=0.9\columnwidth,height=5cm]{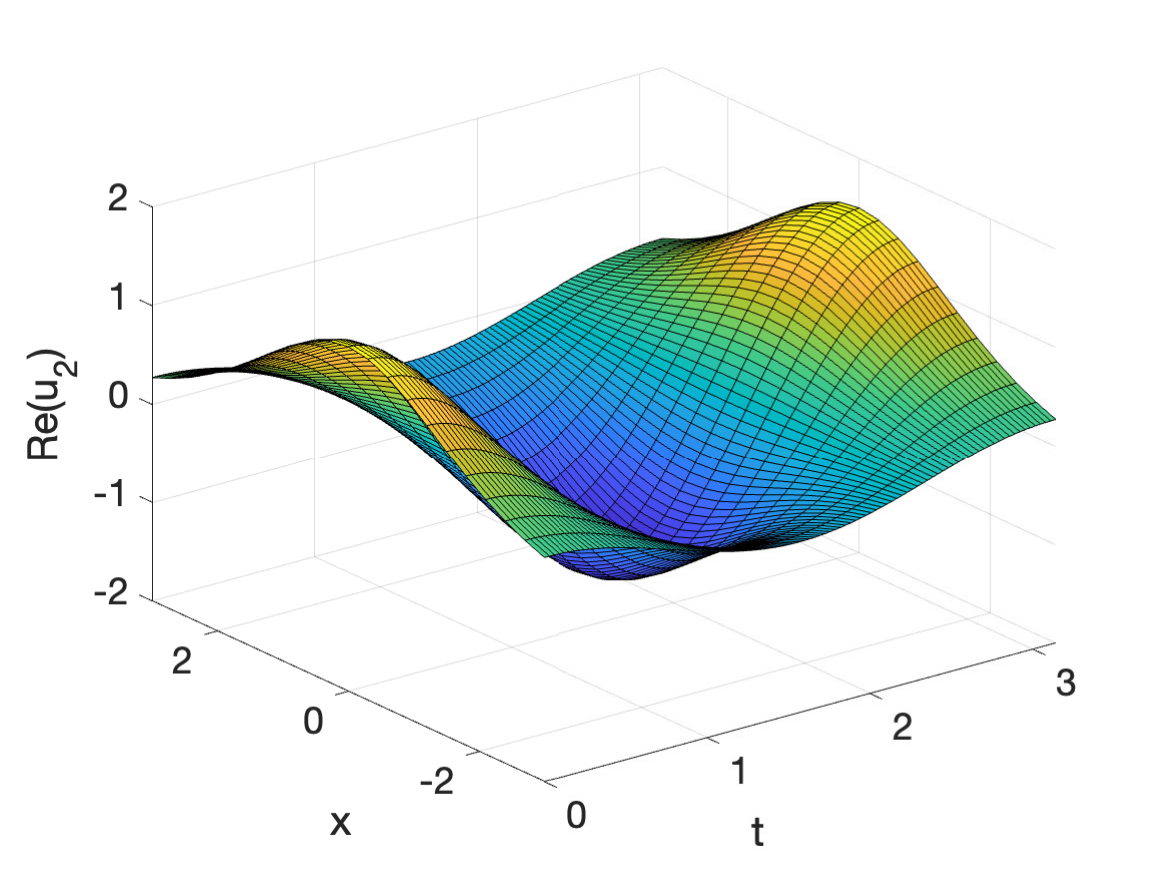}
\par\end{centering}
\centering{}\caption{Plots of the solution to the double simulton using the Fourier spectral
derivative (FSD) method with non-periodic boundaries for the {[}D-D;N-N{]}
boundary condition for the fields $u_{1}$ and $u_{2}$ with $\Delta x=0.3$
and $\Delta t=\pi/2000$.\label{fig:Double_simulton}}
\end{figure}

\begin{table}[!h]
\begin{centering}
\begin{tabular}{|c|c|c|c|c|}
\hline 
Boundary & FSD error & Tc (s) & FIP error & Tc (s)\tabularnewline
\hline 
\hline 
D-D;N-N & \textcolor{red}{$\mathbf{4.5\times10^{-4}}$} & \textcolor{red}{$\mathbf{1.6}$} & \textcolor{red}{$\mathbf{4.5\times10^{-4}}$} & \textcolor{red}{$\mathbf{1.2}$}\tabularnewline
\hline 
N-N;D-N & \textcolor{red}{$\mathbf{4.7\times10^{-4}}$} & \textcolor{red}{$\mathbf{2.1}$} & \textcolor{red}{$\mathbf{4.7\times10^{-4}}$} & \textcolor{red}{$\mathbf{1.4}$}\tabularnewline
\hline 
D-N;N-D & \textcolor{red}{$\mathbf{4.4\times10^{-4}}$} & \textcolor{red}{$\mathbf{2.4}$} & \textcolor{red}{$\mathbf{4.4\times10^{-4}}$} & \textcolor{red}{$\mathbf{1.5}$}\tabularnewline
\hline 
N-D;D-D & \textcolor{red}{$\mathbf{4.1\times10^{-4}}$} & \textcolor{red}{$\mathbf{1.9}$} & \textcolor{red}{$\mathbf{4.1\times10^{-4}}$} & \textcolor{red}{$\mathbf{1.3}$}\tabularnewline
\hline 
\end{tabular}
\par\end{centering}
\medskip{}

\begin{centering}
\begin{tabular}{|c|c|c|c|c|}
\hline 
Boundary & Tau error & Tc (s) & Gal error & Tc (s)\tabularnewline
\hline 
\hline 
D-D;N-N & \textcolor{black}{$7.9\times10^{-6}$} & $2.1$ & \textcolor{black}{$2.1\times10^{-2}$} & $5.6$\tabularnewline
\hline 
N-N;D-N & \textcolor{black}{$1.5\times10^{-5}$} & $2.1$ & \textcolor{black}{$2.2\times10^{-2}$} & $5.4$\tabularnewline
\hline 
D-N;N-D & \textcolor{black}{$9.4\times10^{-6}$} & $2.1$ & \textcolor{black}{$1.9\times10^{-2}$} & $5.8$\tabularnewline
\hline 
N-D;D-D & \textcolor{black}{$3.6\times10^{-6}$} & $2.2$ & \textcolor{black}{$1.8\times10^{-2}$} & $5.6$\tabularnewline
\hline 
\end{tabular}
\par\end{centering}
\centering{}\caption{Errors from solving the double simulton using the FSD, FIP, Tau and
Galerkin methods with four different boundary conditions. The number
of spatial steps is 20 and time steps is 2000 with $\Delta t=\pi/2000$,
$\Delta x=0.3$. The Tau results were obtained using a Legendre basis
with dealiasing factor $d=1$. All errors are normalized by the maximum
value $m=3/2$.\label{tab:Double_simulton}}
\end{table}

In this case, due to the choice of parameters leading to a slowly
varying solution in space, the Tau method gave the lowest errors,
with the FSD and FIP methods having errors between the Tau and Galerkin
methods. These results are tabulated in Table \ref{tab:Double_simulton}. 

\subsection{Summary}

With the Fourier interaction picture (FIP), the error reduction is
extremely large if the problem consists only of linear terms, as seen
for the heat equation. Both FIP and FSD methods yield smaller errors
than a polynomial basis for cases with a rapidly varying spatial dependence,
as observed in the examples of the Peregrine solitary wave and breather
problems.

When the solution is slowly varying compared to the spatial discretization
length, the Tau method with a polynomial basis can produce a smaller
error, as observed in the double simulton problem. The Galerkin method
produced much larger errors than the other spectral methods in all
cases.

The FIP method produces similar errors to the FSD method, with a lower
computational runtime. An exception to this was the very rapidly varying
Peregrine solution, where the use of Neumann boundary conditions resulted
in a higher FIP errors than the FSD method, although still much lower
than the other two polynomial spectral methods for these step-sizes.

It is also possible to perform these comparisons using a finite difference
(FD) method. To demonstrate this, we choose a simple finite difference
scheme, the central differencing method. The corresponding results
are presented in Appendix E. Based on the results, the finite difference
(FD) method always produces larger errors than the spectral methods
and is unstable at large time step-sizes.

\subsection{Detailed error and time scaling comparisons\label{subsec:Detailed-error-and_time-peregrine}}

\begin{figure}
\begin{centering}
\includegraphics[width=0.75\columnwidth]{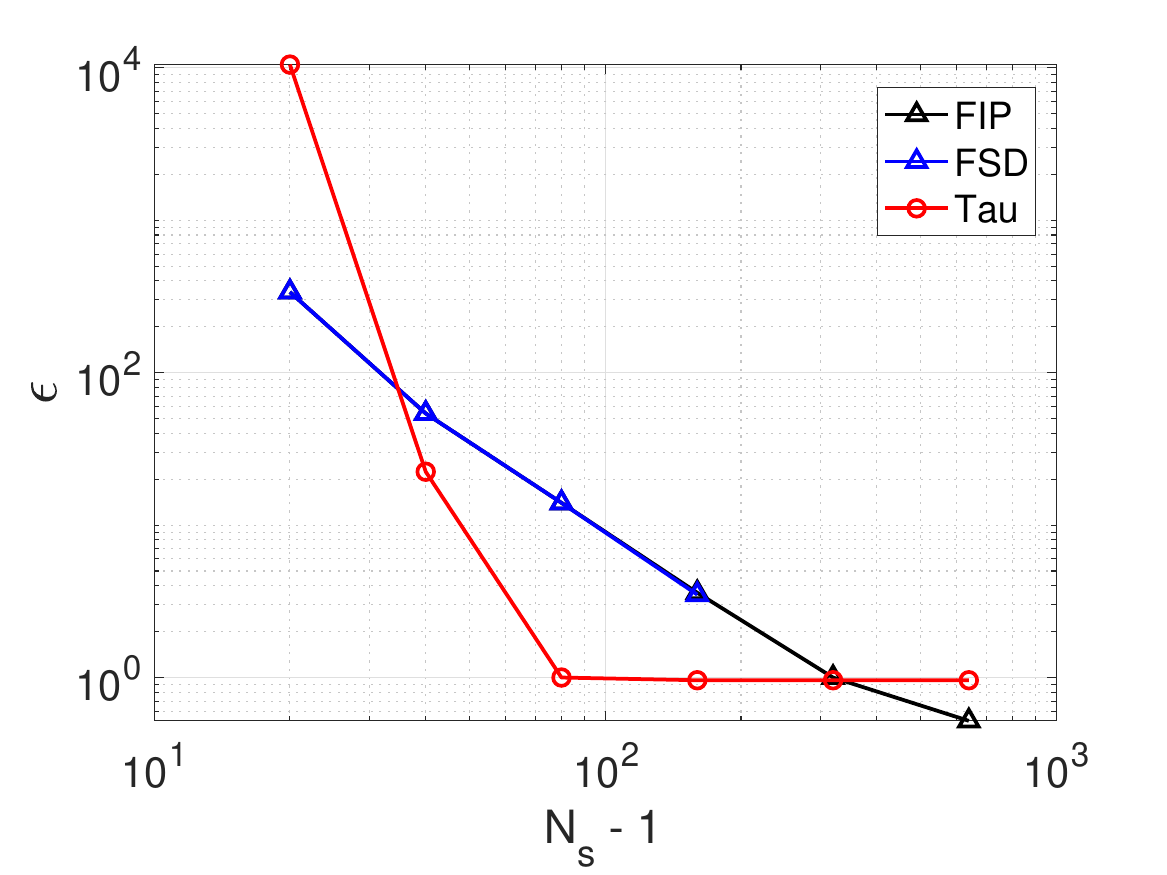}
\par\end{centering}
\begin{centering}
\includegraphics[width=0.75\columnwidth]{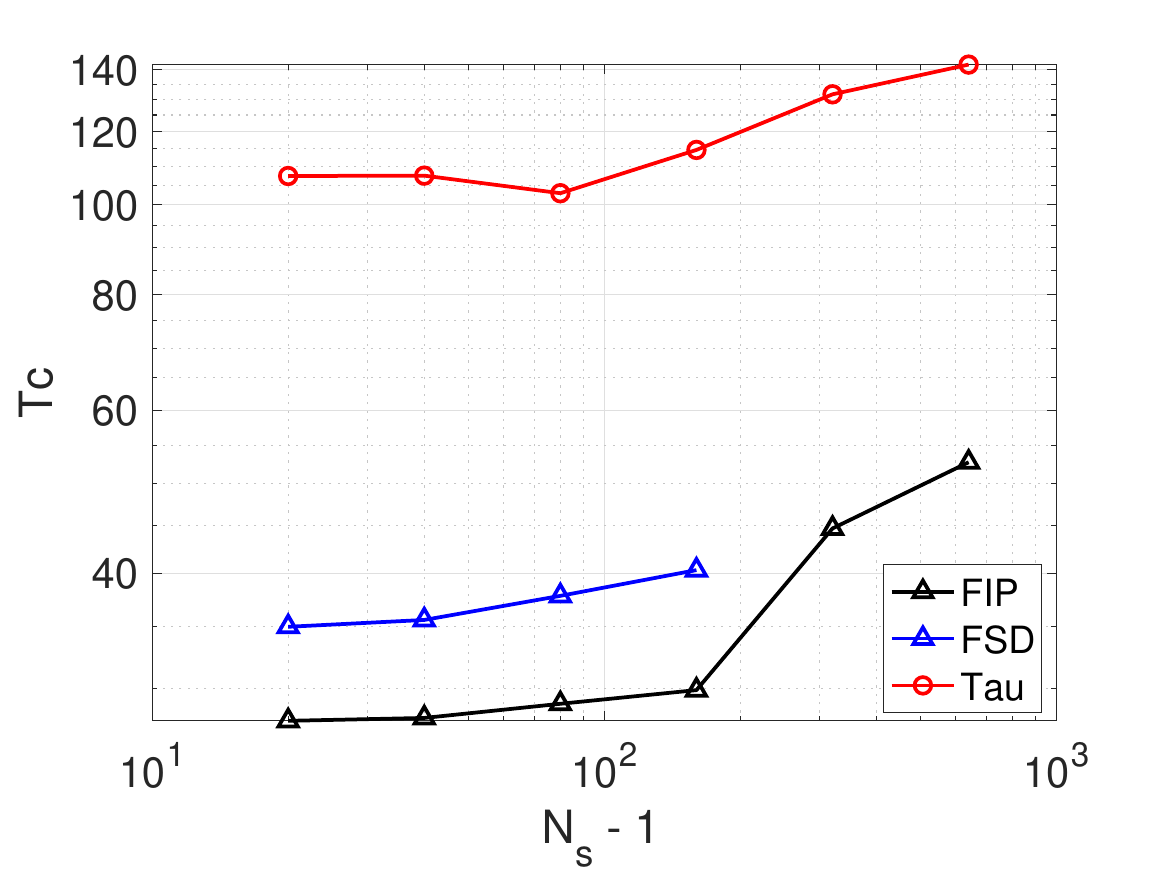}
\par\end{centering}
\caption{The RMS error ($\epsilon$) and time (Tc) scaling comparison for the
Peregrine problem using FIP, FSD and Tau methods with $\triangle t=1.25\times10^{-4}$,
for $\triangle x=0.2,0.1,\ldots0.00625$. The Fourier method is always
faster, but the Tau polynomial method can give lower errors provided
the solution is relatively smooth and the time-step is less than the
von Neumann stability boundary. \label{fig:errorandtimescaling}}
\end{figure}

We now analyze the scaling properties with decreasing spatial grid
step-size, focusing on the two Fourier algorithms and the Tau method.
The spatial step-size error and computation time for the Peregrine
problem are analyzed. The time-step in Table \ref{tab:PeregrineScaling40}
is $20$ times shorter than in the main text examples, with $\triangle t=2.5\times10^{-4}$,
while in Table \ref{tab:PeregrineScaling80} is $40$ times shorter,
with $\triangle t=1.25\times10^{-4}$, in order to demonstrate the
effect of reducing the time-step on the spatial grid scaling. The
equation is the nonlinear Schrödinger equation with initial conditions
that correspond to the Peregrine isolated peak:
\begin{align}
\frac{\partial u}{\partial t} & =i\cdot\left(u\cdot\left|u\right|^{2}+\frac{1}{2}\frac{\partial^{2}u}{\partial x^{2}}\right)\,.
\end{align}
The boundary condition is set to be Dirichlet-Dirichlet. We compute
the time evolution of the field intensity $|u|^{2}$ for the D-D boundary
condition, and obtain an RMS error as tabulated in Tables \ref{tab:PeregrineScaling40}
and \ref{tab:PeregrineScaling80}. The RMS error and time scaling
comparison using the FIP, FSD and Tau methods are presented in Fig.
(\ref{fig:errorandtimescaling}). We note that time-dependent boundary
conditions can, in general, lead to order reduction and the formation
of spurious boundary layers in time-stepping schemes for initial--boundary
value problems, as discussed in \citep{Rosales2024Spatial}. These
issues are known to occur for multi-stage Runge-Kutta schemes, and
are not observed when the midpoint algorithm is used. 
\begin{table}[h]
\centering{}%
\begin{tabular}{|c|c|c|c|c|c|c|}
\hline 
$N_{s}-1$ & FIP: $10^{6}\epsilon$ & Tc (s) & FSD: $10^{6}\epsilon$ & Tc (s) & Tau: $10^{6}\epsilon$ & Tc (s)\tabularnewline
\hline 
\hline 
$20$ & \textbf{\textcolor{red}{$\mathbf{338}$}} & \textcolor{red}{$\bm{13.7}$} & \textcolor{red}{$\mathbf{338}$} & \textcolor{red}{$\mathbf{17.7}$} & $1.05\times10^{4}$ & $53.6$\tabularnewline
\hline 
$40$ & \textbf{\textcolor{red}{$\mathbf{55}$}} & \textcolor{red}{$\mathbf{13.9}$} & \textcolor{red}{$\mathbf{55}$} & \textcolor{red}{$\mathbf{17.9}$} & $22.5$ & $53.7$\tabularnewline
\hline 
$80$ & \textbf{\textcolor{red}{$\mathbf{14}$}} & \textcolor{red}{$\mathbf{14.3}$} & \textcolor{red}{$\mathbf{14}$} & \textcolor{red}{$\mathbf{18.2}$} & $2.78$ & $53.9$\tabularnewline
\hline 
$160$ & \textbf{\textcolor{red}{$\mathbf{4.1}$}} & \textcolor{red}{$\mathbf{15.3}$} & \textcolor{red}{$\mathbf{3.6}$} & \textcolor{red}{$\mathbf{20.1}$} & $2.78$ & $55.9$\tabularnewline
\hline 
$320$ & \textbf{\textcolor{red}{$\mathbf{1.7}$}} & \textcolor{red}{$\mathbf{22.3}$} & \textcolor{red}{$\bm{\mathbf{\infty}}$} & \textcolor{red}{$\mathbf{37.2}$} & $2.78$ & $63.2$\tabularnewline
\hline 
$640$ & \textbf{\textcolor{red}{$\mathbf{1.3}$}} & \textcolor{red}{$\mathbf{27.2}$} & \textcolor{red}{$\bm{\mathbf{\infty}}$} & \textcolor{red}{$\mathbf{48.3}$} & $2.78$ & $68.8$\tabularnewline
\hline 
\end{tabular}\caption{The number of spatial steps, the corresponding relative RMS errors
$\epsilon$ scaled by $10^{6}$, and computation times $Tc$ for the
Peregrine problem in subsection \ref{subsec:Peregrine-solitary-wave}
with $\triangle t=2.5\times10^{-4}$, for $\triangle x=0.2,0.1,\ldots0.00625$.
Algorithms are FIP, FSD and Tau. The Tau results were obtained using
a Chebyshev basis with dealiasing factor $d=2$.\label{tab:PeregrineScaling40}}
\end{table}

\begin{table}[h]
\centering{}%
\begin{tabular}{|c|c|c|c|c|c|c|}
\hline 
$N_{s}-1$ & FIP: $10^{6}\epsilon$ & Tc (s) & FSD: $10^{6}\epsilon$ & Tc (s) & Tau: $10^{6}\epsilon$ & Tc (s)\tabularnewline
\hline 
\hline 
$20$ & \textbf{\textcolor{red}{$\mathbf{338}$}} & \textcolor{red}{$\bm{27.7}$} & \textcolor{red}{$\mathbf{338}$} & \textcolor{red}{$\mathbf{35}$} & $1.05\times10^{4}$ & $107.4$\tabularnewline
\hline 
$40$ & \textbf{\textcolor{red}{$\mathbf{54}$}} & \textcolor{red}{$\mathbf{27.9}$} & \textcolor{red}{$\mathbf{54}$} & \textcolor{red}{$\mathbf{35.6}$} & $22.4$ & $107.5$\tabularnewline
\hline 
$80$ & \textbf{\textcolor{red}{$\mathbf{14}$}} & \textcolor{red}{$\mathbf{28.9}$} & \textcolor{red}{$\mathbf{14}$} & \textcolor{red}{$\mathbf{37.8}$} & $1.00$ & $102.9$\tabularnewline
\hline 
$160$ & \textbf{\textcolor{red}{$\mathbf{3.6}$}} & \textcolor{red}{$\mathbf{29.9}$} & \textcolor{red}{$\mathbf{3.5}$} & \textcolor{red}{$\mathbf{40.3}$} & $0.96$ & $114.6$\tabularnewline
\hline 
$320$ & \textbf{\textcolor{red}{$\mathbf{1.0}$}} & \textcolor{red}{$\mathbf{44.7}$} & \textcolor{red}{$\bm{\mathbf{\infty}}$} & \textcolor{red}{$\mathbf{69.8}$} & $0.96$ & $131.6$\tabularnewline
\hline 
$640$ & \textbf{\textcolor{red}{$\mathbf{0.52}$}} & \textcolor{red}{$\mathbf{52.7}$} & \textcolor{red}{$\bm{\mathbf{\infty}}$} & \textcolor{red}{$\bm{81.9}$} & $0.96$ & $141.7$\tabularnewline
\hline 
\end{tabular}\caption{The number of spatial steps, the corresponding relative RMS errors
$\epsilon$ scaled by $10^{6}$, and computation times $Tc$ for the
Peregrine problem in subsection \ref{subsec:Peregrine-solitary-wave}
with $\triangle t=1.25\times10^{-4}$, for $\triangle x=0.2,0.1,\ldots0.00625$.
Algorithms are FIP, FSD and Tau. The Tau results were obtained using
a Chebyshev basis with dealiasing factor $d=2$.\label{tab:PeregrineScaling80}}
\end{table}

In Table \ref{tab:PeregrineScaling80}, the error scaling of both
Fourier methods approximately quadratic, $\epsilon\propto\Delta x^{2}$,
for space-steps $\Delta x\ge0.025$, or $N_{s}\le160$. In this regime
errors are independent of time-step. For smaller space-steps, the
errors are partly due to time discretization. In this regime the reduction
in FIP error is less than quadratic, becoming linear at $N_{s}\approx160$,
and eventually increasing if the space-step is reduced further.

In this regime the maximum Laplacian eigenvalue is so large that the
transformed ODE equations require extremely small time-steps, if relative
error levels below $10^{-6}$ are required. However, for the FSD method
at smaller space-steps, the errors increase rapidly to the point where
the algorithm is completely divergent, giving infinite errors. These
effects are seen in both tables.

The von Neumann stability boundary for an FSD scheme for the linear
heat equation is $\Delta t\le\Delta x^{2}$, which corresponds to
$\Delta x\ge0.011$ in the upper table and $\Delta x\ge0.079$ in
the lower table. Thus, the FSD method has a divergence at space-steps
slightly larger than the linear criterion. There is no instability
boundary for the FIP method.

An unexpected feature of these results is the sub-linear behavior
of the space-step dependence, caused partly by the vector nature of
current CPU arithmetic units and partly by overheads in the software.
For both Fourier methods, a large part of the floating-point arithmetic
can be carried out in parallel. Hence an $N\log N$ complexity scaling
is not seen for small grid sizes, with similar results for the Tau
method.

For both methods, only for $N_{S}>160$ is there a change in timing
due to the number of spatial points. Even for $N_{S}=$ 640 , the
Fourier scaling is still sub-linear, with $T\propto\sqrt{N_{s}}$,
although for extremely large grids one expects asymptotic $N_{S}\log N_{S}$
behavior. However, since there is no parallelization of the time-stepping
part of the algorithm, the time taken in all cases is strictly linear
in the number of time-steps.

\subsection{Time-step convergence and complexity order}

Convergence depends on both time-step and spatial discretization.
To test the convergence while accounting for this coupling between
the time and space steps, one must minimize the error at a fixed resource.
For this, the product of time and space steps is held fixed while
the relative number of time steps is changed. This procedure is then
repeated for different fixed resources. This can be used to obtain
the complexity order, discussed in \citep{teh2025complexity}, on
how the total error scales with total resource use. 

For independent resources, there is a known result for the complexity
order $c,$in terms of individual convergence orders $n_{i}$ \citep{teh2025complexity},
namely: 
\begin{equation}
c^{-1}=\sum_{i}n_{i}^{-1}
\end{equation}
This shows that the complexity order is bounded by the lowest order
for any contributing resource. Since the FIP method allows the use
of other higher-order time-stepping schemes, such as the fourth-order
Runge-Kutta (RK4) algorithm \citep{Holmes2007,caradoc2000vortex,Butcher1987thenumerical}
we present results obtained using the FIP method with both MP and
RK4 algorithms, and compare them with results obtained using the Tau
method with RK222 and RK443 \citep{Ascher1997implicit} schemes.

\begin{figure}
\begin{centering}
\includegraphics[width=0.75\columnwidth]{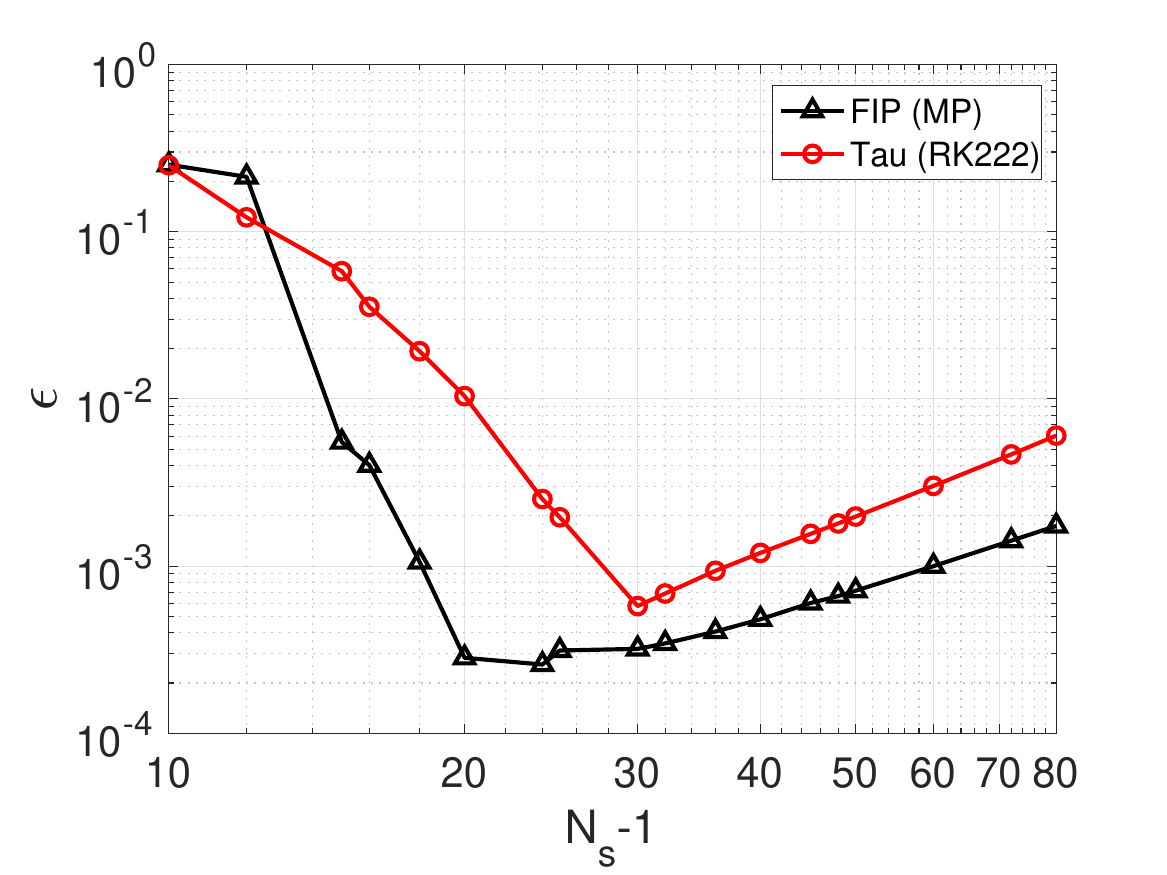}
\par\end{centering}
\begin{centering}
\includegraphics[width=0.75\columnwidth]{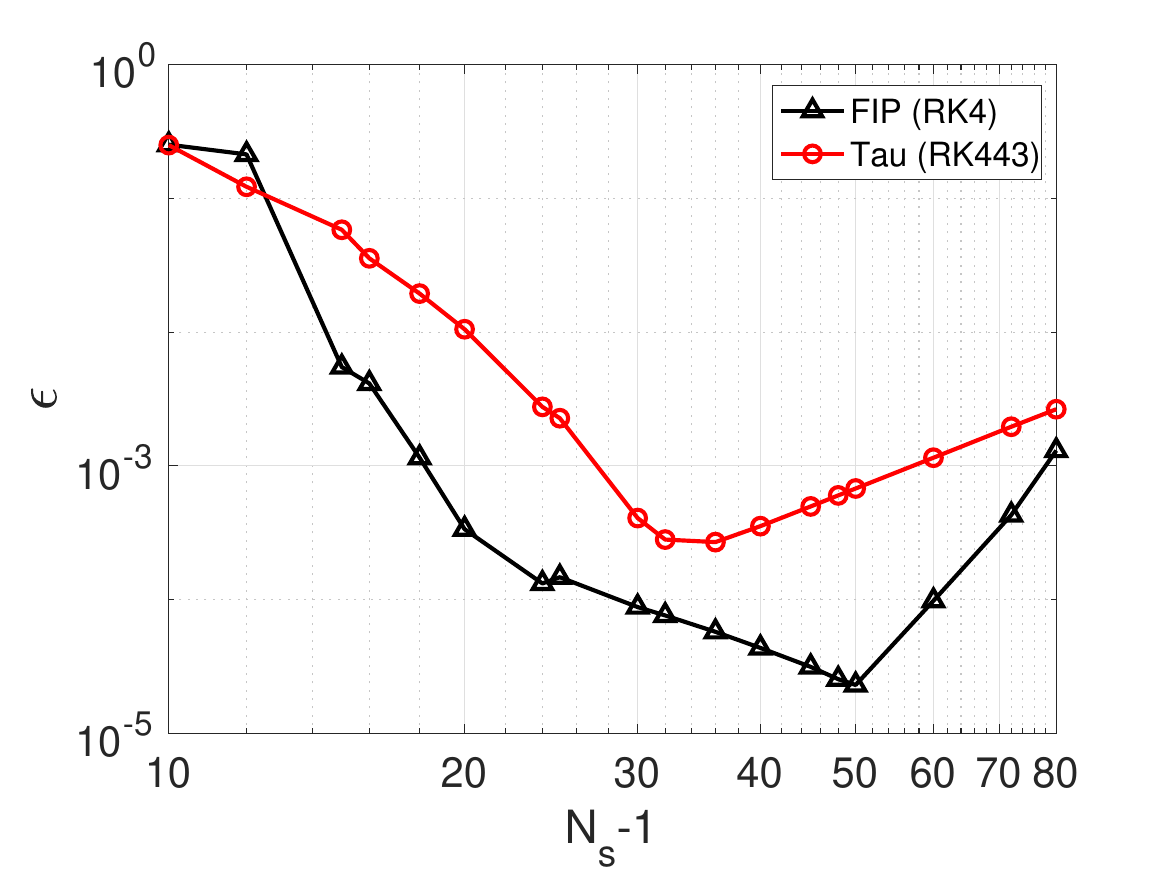}
\par\end{centering}
\caption{RMS error ($\epsilon$) comparisons for the Peregrine problem using
FIP and Tau methods for a fixed resource of $N_{R}=72000$. The upper
graph compares FIP (MP) and Tau (RK222), the lower graph compares
FIP (RK4) and Tau (RK443). A clear local minimum is observed given
a fixed resource $N_{R}$. The Tau results used a Chebyshev basis
with dealiasing factor $d=2$. \label{fig:varying_resources_N=00003D72000}}
\end{figure}

Here we analyze the RMS error $\epsilon$ for a fixed resource $N_{R}=(N_{s}-1)(N_{t}-1)$
for the Peregrine problem, where $N_{s}-1$ is the number of space
steps and $N_{t}-1$ is the number of time steps. The RMS error comparison
for a fixed resource of $N_{R}=72000$ and $N_{R}=100800$, obtained
using the FIP and Tau methods with different time-stepping schemes
is presented in Fig. (\ref{fig:varying_resources_N=00003D72000}).

Based on Fig (\ref{fig:varying_resources_N=00003D72000}) , a local
minimum for the RMS error $\epsilon$ is evident for a fixed resource
$N_{R}$ for all methods. The minimum RMS error ($\epsilon_{min}$)
obtained by varying the fixed resource $N_{R}$ using the FIP and
Tau methods with different time-stepping schemes, is presented in
Table \ref{tab:varying_resource_minimum}.

\begin{table}[H]
\begin{centering}
\begin{tabular}{|c|c|c|c|c|}
\hline 
$N_{R}/1k$ & \begin{cellvarwidth}[t]
\centering
FIP-MP\\
 $10^{4}\epsilon_{min}$
\end{cellvarwidth} & \begin{cellvarwidth}[t]
\centering
FIP-RK4\\
 $10^{4}\epsilon_{min}$
\end{cellvarwidth} & \begin{cellvarwidth}[t]
\centering
Tau-Rk222\\
 $10^{4}\epsilon_{min}$
\end{cellvarwidth} & \begin{cellvarwidth}[t]
\centering
Tau-RK443\\
 $10^{4}\epsilon_{min}$
\end{cellvarwidth}\tabularnewline
\hline 
\hline 
7.2 & \textcolor{red}{$\mathbf{46.1}$} & \textcolor{red}{$\mathbf{2.39}$} & $381$ & $271$\tabularnewline
\hline 
18 & \textcolor{red}{$\mathbf{8.19}$} & \textcolor{red}{$\mathbf{0.48}$} & $81.1$ & $52.7$\tabularnewline
\hline 
25.2 & \textcolor{red}{$\mathbf{5.97}$} & \textcolor{red}{$\mathbf{0.38}$} & $33.6$ & $27.7$\tabularnewline
\hline 
36 & \textcolor{red}{$\mathbf{3.75}$} & \textcolor{red}{$\mathbf{0.30}$} & $16.7$ & $12.2$\tabularnewline
\hline 
50.4 & \textcolor{red}{$\mathbf{2.73}$} & \textcolor{red}{$\mathbf{0.26}$} & $10.3$ & $5.67$\tabularnewline
\hline 
72 & \textcolor{red}{$\mathbf{2.58}$} & \textcolor{red}{$\mathbf{0.23}$} & $5.77$ & $2.69$\tabularnewline
\hline 
100.8 & \textcolor{red}{$\mathbf{1.92}$} & \textcolor{red}{$\mathbf{0.22}$} & $3.40$ & $1.04$\tabularnewline
\hline 
\end{tabular}
\par\end{centering}
\caption{Minimum RMS error $\epsilon_{min}$ scaled by $10^{4}$, for varying
the fixed resource $N_{R}$ from $7.2\times10^{3}$ to $100.8\times10^{3}$
for the Peregrine problem. Results are obtained for the FIP (MP, RK4)
and Tau (RK222, RK443) methods. The Tau results used a Chebyshev basis
with dealiasing factor $d=2$. \label{tab:varying_resource_minimum}}
\end{table}

\begin{figure}[H]
\begin{centering}
\includegraphics[width=0.75\columnwidth]{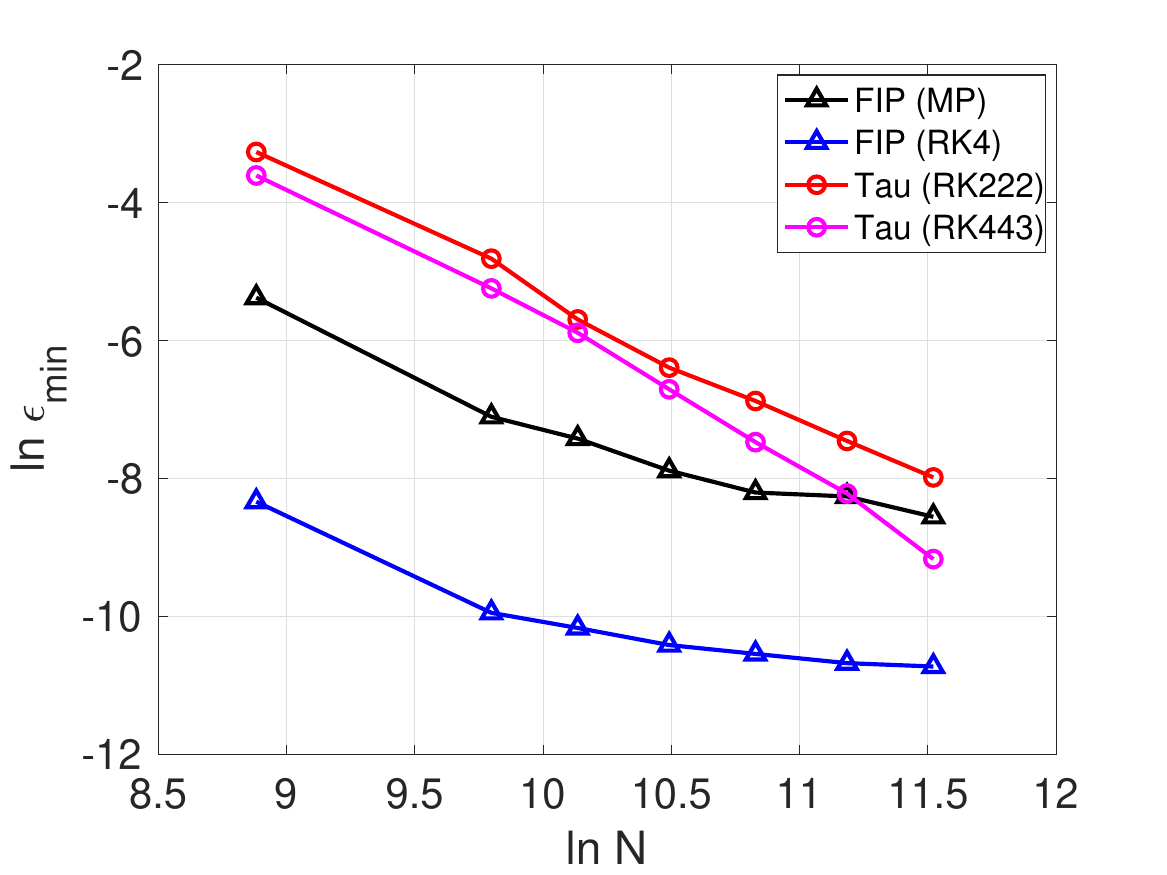}
\par\end{centering}
\caption{Comparison of the minimum RMS error ($\epsilon_{min}$), for fixed
resource $N_{R}$ for the Peregrine problem, using the FIP (MP, RK4)
and Tau (RK222, RK443) methods. The Fourier method generates lower
$\epsilon_{min}$, compared to the Tau polynomial method when using
the same temporal orders. \label{fig:minimum_vs_resource}}
\end{figure}

From the results in Fig. (\ref{fig:minimum_vs_resource}), the Fourier
method (MP, RK4) generates lower errors with the same temporal order
methods in all cases. While the polynomial spectral method may give
lower errors at very large resource use, this requires extremely large
resources to reach the asymptotic region. 

\section{stochastic partial differential equations \label{sec:stochastic-heat}}

Finally, we consider stochastic partial differential equations (SPDE)
with a non-periodic boundary condition and show that the Fourier algorithm
is especially useful in solving this type of SPDE. We solve the (1+1)-dimensional
stochastic heat equation, an SPDE with an analytic solution, using
both the FIP and the Tau methods for comparison. In addition, we analyze
the error and time scaling properties of this SPDE, similar to the
comparisons discussed in subsection \ref{subsec:Detailed-error-and_time-peregrine}.
We also solve a (1+2)-dimensional stochastic heat equation with correlated
noise and an analytic solution to demonstrate that the Fourier algorithm
extends naturally to more complex higher-dimensional problems.

\subsection{\textcolor{black}{Stochastic heat equation}\textcolor{blue}{\label{sec:Non-periodic-boundary-condition}}}

We treat an exactly soluble stochastic heat equation \citep{Bertini1995,Alos1999,Swanson2007}
given by
\begin{align}
\frac{\partial u\left(t,x\right)}{\partial t} & =\frac{1}{2}\frac{\partial^{2}u\left(t,x\right)}{\partial x^{2}}+\eta\left(t,x\right)\,,\label{eq:SPDE}
\end{align}
where the noise term $\eta(t,x)$ is delta correlated in space and
time, with the noise correlation $\langle\eta_{i}(t,x)\eta_{j}(t',x')\rangle=\delta(t-t')\delta(x-x')$.
We impose zero boundary conditions, so that, for $k=\pi/L$,
\begin{equation}
u(t,x)=\sum_{n\ge1}u_{n}\left(t\right)\sin\left(nkx\right)=\sum_{n>1}u_{n}\left(t\right)\sin\left(n\pi x/L\right)
\end{equation}
have vanishing boundaries at $x=0,\,L$. By applying the discrete
sine transform to the SPDE, we get the equation
\begin{align}
\frac{\partial u_{n}\left(t\right)}{\partial t} & =-K_{n}u_{n}\left(t\right)+\eta_{n}\left(t\right)\,,\label{eq:transform_eqn}
\end{align}
where $K_{n}=n^{2}\pi^{2}/(2L^{2})$ and 
\begin{equation}
\eta_{n}\left(t\right)=\frac{2}{L}\int_{0}^{L}\eta\left(t,x\right)\sin\left(n\pi x/L\right)dx\,.
\end{equation}
The transformed SPDE Eq. (\ref{eq:transform_eqn}) has a solution:
\begin{equation}
u_{n}\left(t\right)=e^{-K_{n}t}\left(u_{n}\left(0\right)+\int_{0}^{t}e^{K_{n}\tau}\eta_{n}\left(\tau\right)d\tau\right)\,.
\end{equation}
Therefore,
\begin{align}
u\left(t,x\right) & =\sum_{n\ge1}\sin\left(n\pi x/L\right)e^{-K_{n}t}\nonumber \\
 & \;\times\left(u_{n}\left(0\right)+\int_{0}^{t}e^{K_{n}\tau}\eta_{n}\left(\tau\right)d\tau\right)\,.
\end{align}
With this expression for $u(t,x)$ and assuming the initial conditions
$u_{n}(0)=0$, an observable 
\begin{align}
J(t) & \equiv\intop_{0}^{L}\langle\left[u(t,x)\right]{}^{2}\rangle\,dx\label{eq:Integrated-mean-square}
\end{align}
 can be shown to have an analytical solution
\begin{align}
J\left(t\right) & =\sum_{n\ge1}\frac{L^{2}}{n^{2}\pi^{2}}\left(1-e^{-\frac{n^{2}\pi^{2}}{L^{2}}t}\right)\,.\label{eq:J(t)-1D}
\end{align}

We ran numerical simulations using the FIP and Tau methods, choosing
the time to range from $0$ to $1$, while the spatial interval was
from $0$ to $L=5$. The analytical solution was summed numerically
up to $n=10^{7}$ to guarantee convergence. Numerical results are
tabulated in Table \ref{tab:numerics_SPDE}.

\begin{figure}[H]
\centering{}\includegraphics[width=0.75\columnwidth]{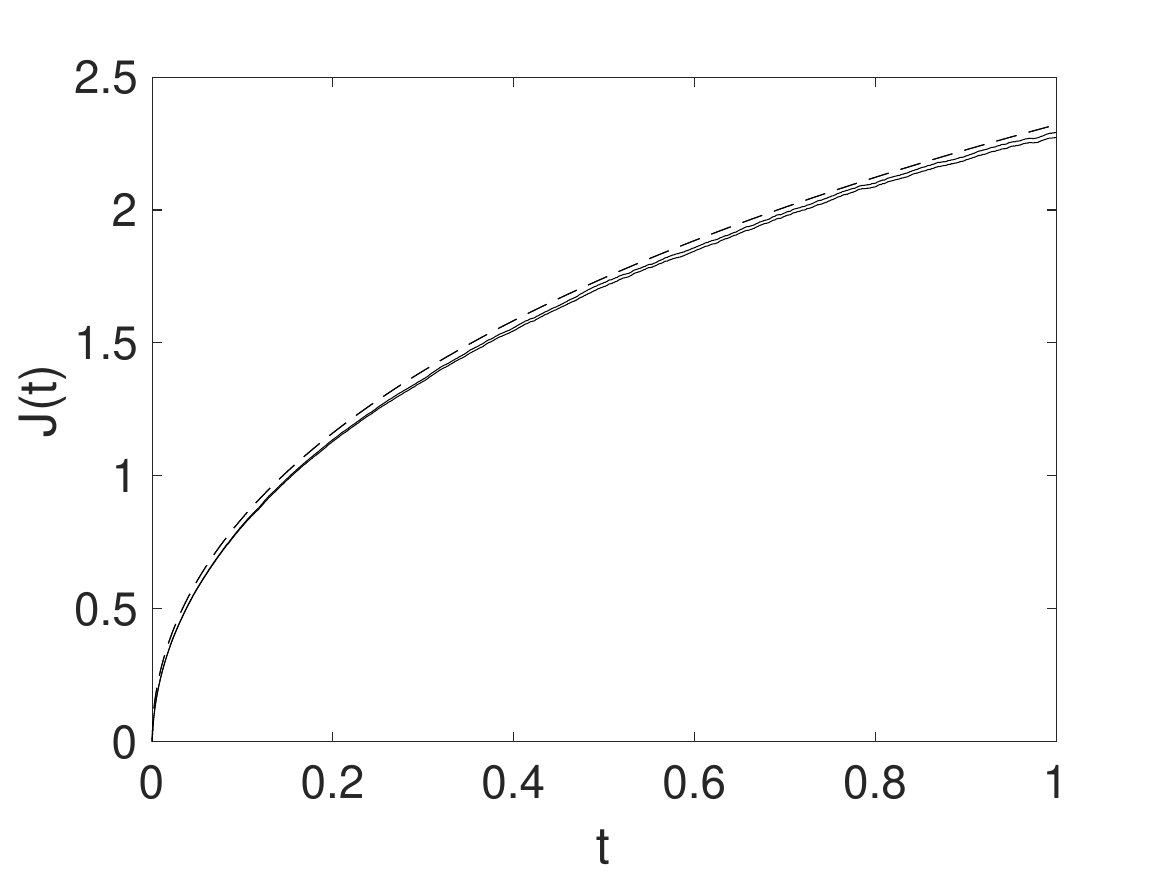}\caption{Plot of the integrated mean square solution $J(t)$ of the stochastic
heat equation, computed using the FIP method for the D-D boundary
condition with $\Delta x=0.05$, $\Delta t=1/1000$ and $2\times10^{4}$
samples. The solid lines indicate upper and lower sampling error bounds.
The dashed line indicates the analytical solution in Eq. (\ref{eq:J(t)-1D})
summed up to $n=10^{7}$.\label{fig:Fig_SPDE}}
\end{figure}

\begin{table}[H]
\begin{centering}
\begin{tabular}{|c|c|c|c|c|}
\hline 
Boundary & FIP error & Tc (s) & Tau error & Tc(s)\tabularnewline
\hline 
\hline 
D-D & \textcolor{red}{$\mathbf{1.37\times10^{-2}}$} & \textcolor{red}{$\mathbf{469.7}$} & $2.92\times10^{-2}$ & $2447.2$\tabularnewline
\hline 
\end{tabular}
\par\end{centering}
\centering{}\caption{Errors and timing for solving the SPDE in equation \ref{eq:SPDE}
using the FIP (in red) and Tau methods, showing lower errors and faster
timing. The simulation parameters were $\Delta x=0.05$, $\Delta t=1/1000$,
and $2\times10^{4}$ samples. The Tau results were obtained using
a Chebyshev basis with dealiasing factor $d=1$. The errors are dominated
by the space step error. \label{tab:numerics_SPDE}}
\end{table}

The results in Table \ref{tab:numerics_SPDE} show that the FIP method
has better than $50\%$ lower RMS errors and is more than 5 times
faster for solving the stochastic heat equation. The much larger speed
improvement could be from internal software differences. Unlike the
deterministic equations treated above, using smaller space-steps does
not give a relatively smoother solution, due to the noise terms. A
detailed treatment of the sampling error, time and space-step error
scaling for this problem using the FIP method is given elsewhere \citep{teh2025complexity}.

\subsubsection{Error and time scaling comparisons}

Next we analyze the space-step error scaling and computation time
for the (1+1)-dimensional stochastic heat equation. The time-step
size is fixed at $\Delta t=1\times10^{-3}$, the sampling size at
$2\times10^{4}$, and the number of space-steps is varied from 10
to 80, giving a time-step error of $10^{-4}$ to $10^{-6}$ and a
sampling error of order $10^{-3}$. This allows the space-step error
scaling properties to be analyzed, as in these cases the error is
dominated by the space-step error.

\begin{figure}[H]
\begin{centering}
\textcolor{blue}{\includegraphics[width=0.75\columnwidth]{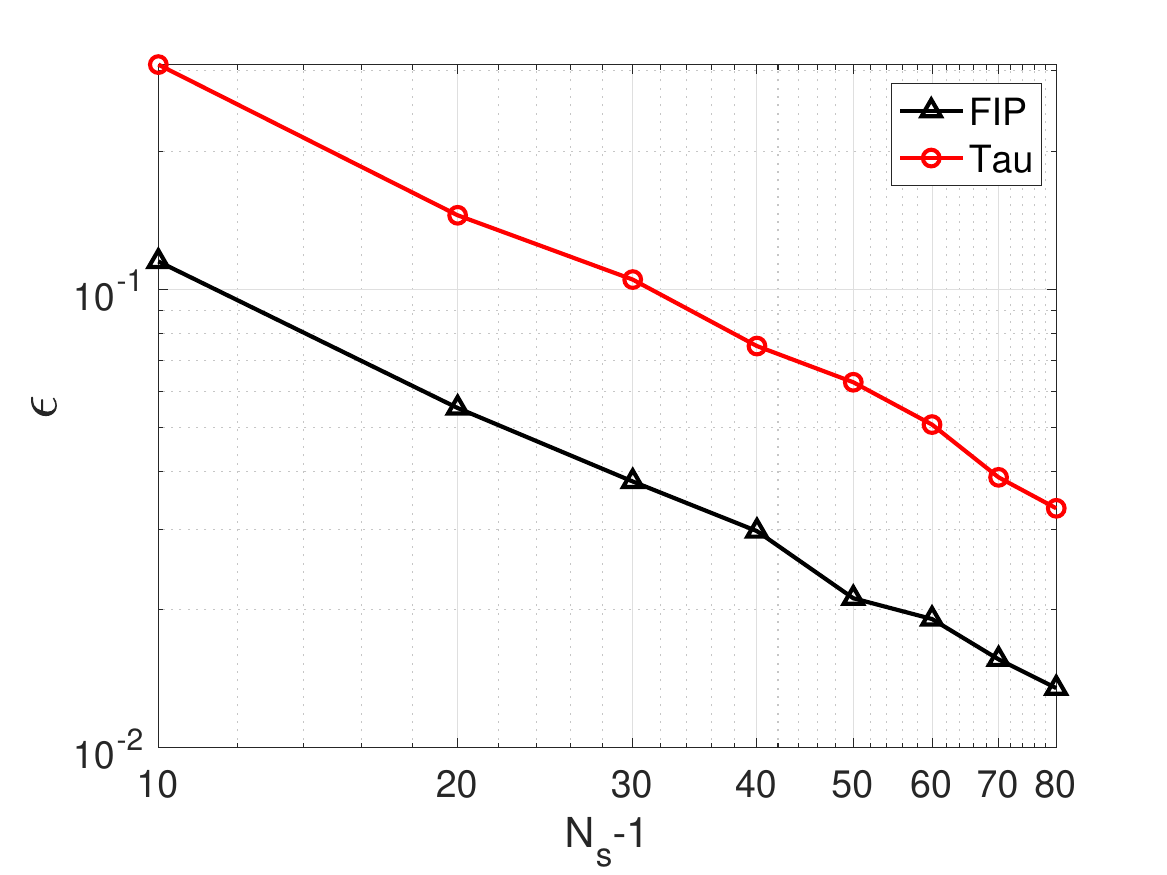}}
\par\end{centering}
\begin{centering}
\textcolor{blue}{\includegraphics[width=0.75\columnwidth]{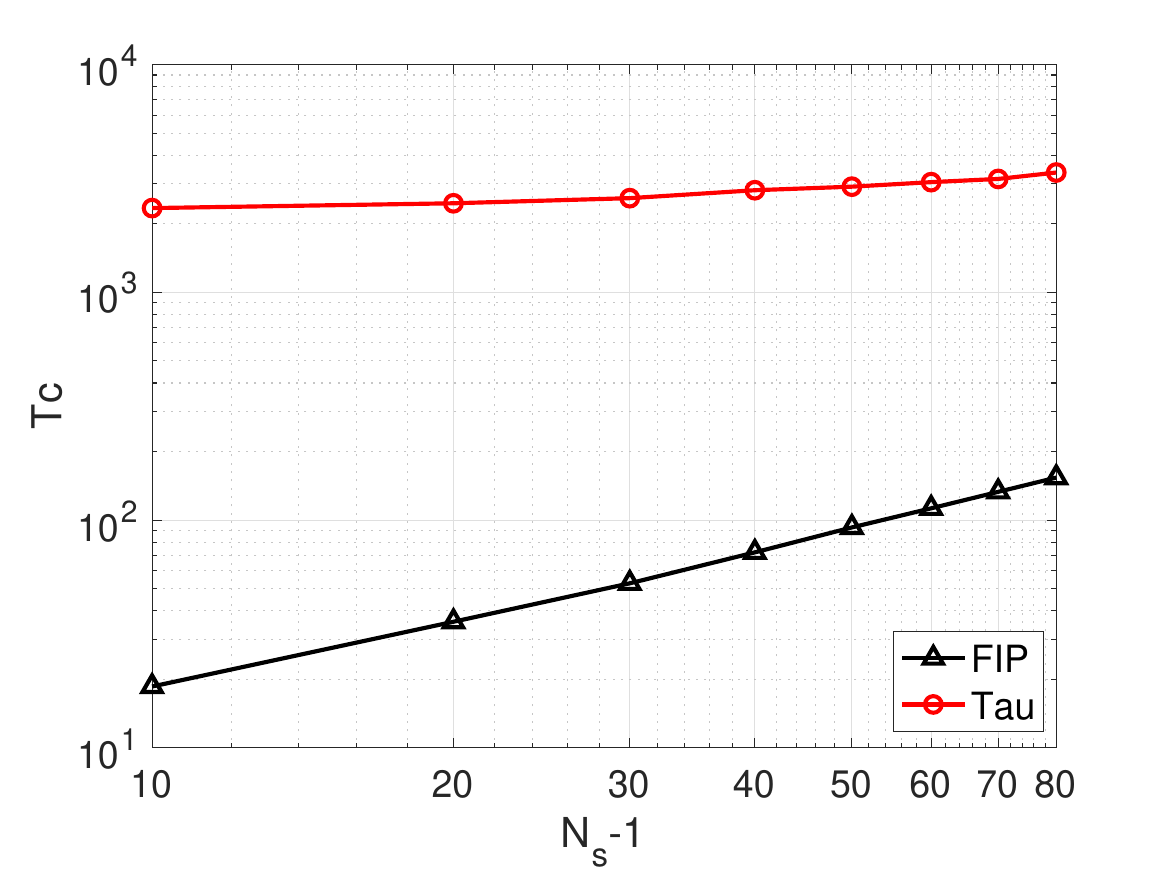}}
\par\end{centering}
\caption{The RMS error ($\epsilon$) and time (Tc) scaling comparison for the
stochastic heat equation using the FIP and Tau methods with $\triangle t=1\times10^{-3}$,
for $\triangle x=0.4,0.2,\ldots0.05$. The Fourier method is faster
and gives lower errors compared to the Tau polynomial method. \label{fig:SPDE_scaling_10-80}}
\end{figure}

We consider the exactly soluble stochastic heat equation discussed
in section \ref{sec:Non-periodic-boundary-condition}. The boundary
condition is set to be Dirichlet-Dirichlet. We compute the observable
$J(t)$ for the D-D boundary condition, and obtain an RMS error as
tabulated in Table \ref{tab:SPDE_scaling_10-80}. The RMS error and
time scaling comparison using the FIP and Tau methods are presented
in Fig. (\ref{fig:SPDE_scaling_10-80}).

\begin{table}[H]
\begin{centering}
\begin{tabular}{|c|c|c|c|c|}
\hline 
$N_{s}-1$ & FIP: $10^{2}\epsilon$ & Tc (s) & Tau: $10^{2}\epsilon$ & Tc (s)\tabularnewline
\hline 
\hline 
10 & \textcolor{red}{$\mathbf{11.5}$} & \textbf{\textcolor{red}{$\mathbf{18.6}$}} & $30.9$ & $2344$\tabularnewline
\hline 
20 & \textcolor{red}{$\mathbf{5.51}$} & \textcolor{red}{$\mathbf{35.8}$} & $14.5$ & $2460.8$\tabularnewline
\hline 
30 & \textcolor{red}{$\mathbf{3.81}$} & \textcolor{red}{$\mathbf{52.8}$} & $10.5$ & $2591.6$\tabularnewline
\hline 
40 & \textcolor{red}{$\mathbf{2.97}$} & \textcolor{red}{$\mathbf{72.2}$} & $7.52$ & $2808$\tabularnewline
\hline 
50 & \textcolor{red}{$\mathbf{2.12}$} & \textcolor{red}{$\mathbf{92.9}$} & $6.27$ & $2911$\tabularnewline
\hline 
60 & \textcolor{red}{$\mathbf{1.91}$} & \textcolor{red}{$\mathbf{112.7}$} & $5.07$ & $3047$\tabularnewline
\hline 
70 & \textcolor{red}{$\mathbf{1.56}$} & \textcolor{red}{$\mathbf{133.1}$} & $3.89$ & $3147.5$\tabularnewline
\hline 
80 & \textcolor{red}{$\mathbf{1.35}$} & \textcolor{red}{$\mathbf{153.7}$} & $3.33$ & $3357.9$\tabularnewline
\hline 
\end{tabular}
\par\end{centering}
\caption{Spatial steps, RMS errors $\epsilon$ scaled by $10^{2}$, and computation
times $Tc$ for the stochastic heat equation in section \ref{sec:Non-periodic-boundary-condition}
with $\triangle t=1\times10^{-3}$, for $\triangle x=0.4,0.2,\ldots0.05$.
Algorithms are FIP and Tau, obtained using a Chebyshev basis with
dealiasing factor $d=1$. \label{tab:SPDE_scaling_10-80}}
\end{table}

Based on Fig. (\ref{fig:SPDE_scaling_10-80}), the FIP method shows
lower RMS errors than the polynomial method across the entire range,
indicating a much higher accuracy for the same number of space-steps.
 Unlike deterministic cases, there is no relative improvement with
polynomial methods at smaller step-sizes, because the solutions do
not become smoother in this limit. For this stochastic partial differential
equation, the spatially delta-correlated noise gives a convergence
order of about $n=1$ in both cases, but with different prefactors.

\subsection{Stochastic heat equation with filtered noise \label{subsec:(1+1)-d-heta_filtered_noise}}

Next, we consider a case where equally-spaced grid-points is especially
useful, an SPDE whose stochastic noise terms have a finite spatial
correlation distance. Finite correlation lengths are typical of many
noise processes in realistic stochastic spatio-temporal theories \citep{Traulsen2004Generation,Sagules2007spatiotemporal,Busch2003influence,Lam1993spatiotemporal},
including interdisciplinary dynamical cases such as geophysical, sociological
and biophysical models. These can be efficiently treated numerically
using Fourier transforms on an equally-spaced grid that with spectral
filtering. 

To give an exactly soluble case, we solve the (1+1)-dimensional stochastic
heat equation of the form in Eq. (\ref{eq:SPDE}), with filtered noise
where the noise term $\eta(t,x)$ is delta correlated in time, but
contains a finite correlation in space, with the noise satisfying
the moment
\begin{align}
\langle\eta(t,x)\eta(t',x')\rangle & =\delta(t-t')\left(\frac{\kappa}{2}\right)\exp(-\kappa|x-x'|)\,.
\end{align}
Here $\kappa$ is a cutoff that defines the correlation length in
space. In the limit $\kappa\rightarrow\infty$, $\langle\eta(t,x)\eta(t',x')\rangle=\delta(t-t')\delta(x-x')$.
The filtered spatial white noise or colored noise provides a more
physically realistic noise with an exponential correlation function
in space, similar to an Ornstein-Uhlenbeck process in time \citep{Uhlenbeck1930onthetheory,Traulsen2004Generation,Gardiner2004handbook}.

We subject this stochastic heat equation to different types on boundary
conditions: one with vanishing boundaries at $x=0,\,L$ (D-D boundary
condition) and the other with vanishing spatial derivative at the
boundaries $x=0,\,L$ (N-N boundary condition). For these types of
boundary conditions, analytical solutions exist for the observable
$J(t)$ in Eq. (\ref{eq:Integrated-mean-square}).

For the D-D boundary conditions, the solution $u(t,x)$ is expanded
as
\begin{align}
u(t,x) & =\sum_{n\geq1}u_{n}(t)\text{sin}(\frac{n\pi x}{L})\,,
\end{align}
and the quantity $J(t)$ can be shown to be
\begin{align}
J(t) & =\sum_{n\geq1}\frac{1-e^{-2K_{n}t}}{2K_{n}}P_{nn}\,,
\end{align}
where 
\begin{align}
P_{nn} & =\frac{1}{1+\left(n\pi/\kappa L\right)^{2}}+\frac{2\left(n\pi/\kappa L\right)^{2}(1-e^{-\kappa L}(-1)^{n})}{\left(\kappa L\right)\left(1+\left(n\pi/\kappa L\right)^{2}\right)^{2}}\,.
\end{align}

On the other hand, for the N-N boundary condition, the solution $u(t,x)$
is expanded as
\begin{align}
u(t,x) & =\sum_{n=0}u_{n}(t)\text{cos}(nkx)\,.
\end{align}
The exact solution for $J(t)$ is
\begin{align}
J(t) & =t\left(1-\frac{1-e^{-\kappa L}}{\kappa L}\right)+\sum_{n\geq0}\frac{1-e^{-2K_{n}t}}{2K_{n}}Q_{nn}\,,\label{eq:J(t)-NN}
\end{align}
where 
\begin{align}
Q_{nn} & =\frac{1}{1+\left(n\pi/\kappa L\right)^{2}}-\frac{2(1-e^{-\kappa L}(-1)^{n})}{\left(\kappa L\right)\left(1+\left(n\pi/\kappa L\right)^{2}\right)^{2}}\,.
\end{align}

We ran numerical simulations using the FIP method, choosing the time
to range from $0$ to $1$, while the spatial interval was from $0$
to $L=4$ with $\kappa=10$. The analytical solution was summed numerically
up to $n=10^{4}$ to guarantee convergence. Numerical results are
tabulated in Table \ref{tab:numerics_SPDE-DD-NN}.

\begin{figure}[H]
\centering{}\includegraphics[width=0.75\columnwidth]{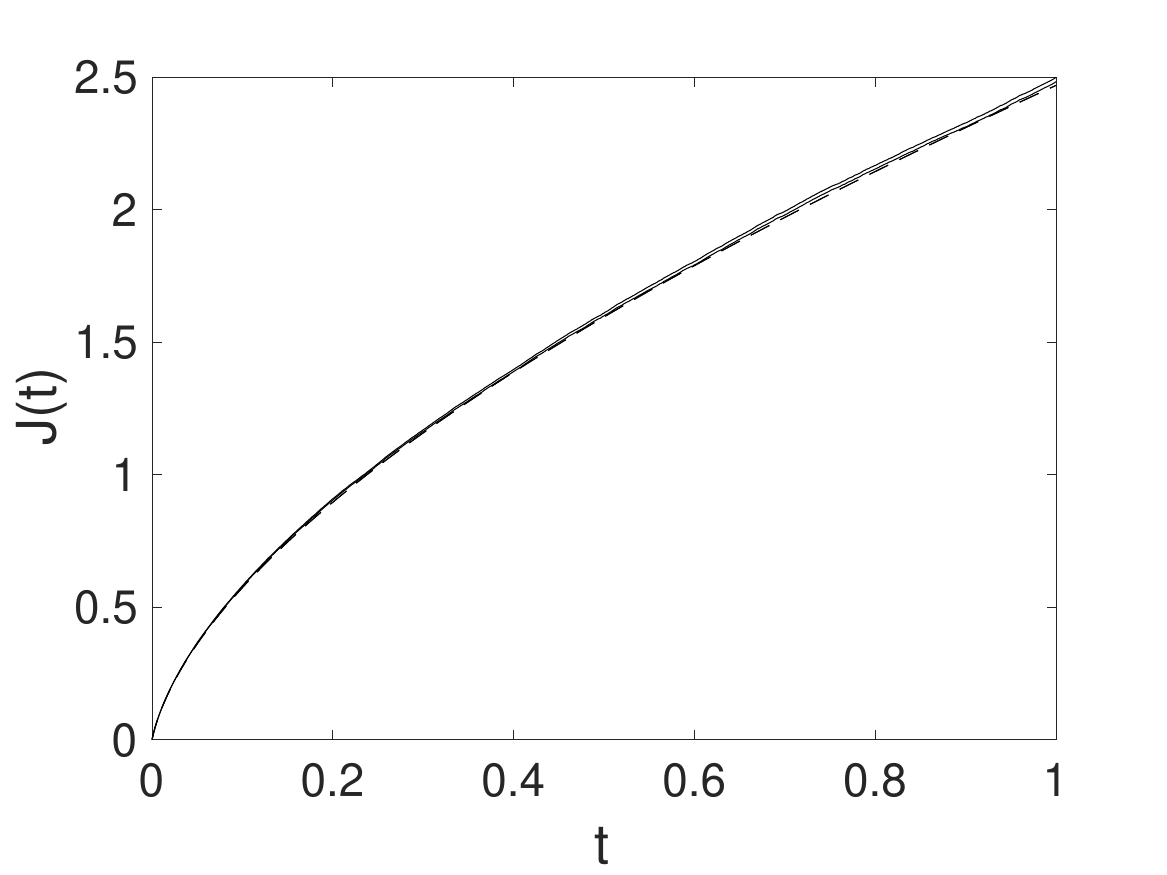}\caption{Plot of the integrated mean square solution $J(t)$ of the stochastic
heat equation with filtered noise, using the FIP method and N-N boundary
condition with $\Delta x=0.04$, $\Delta t=1\times10^{-3}$ and $2\times10^{4}$
samples. The solid lines indicate upper and lower sampling error bounds.
The dashed line is the analytical solution in Eq. (\ref{eq:J(t)-NN})
summed up to $n=10^{4}$.\label{fig:Fig_SPDE-NN}}
\end{figure}

\begin{table}[H]
\begin{centering}
\begin{tabular}{|c|c|c|}
\hline 
Boundary & FIP error & Tc (s)\tabularnewline
\hline 
\hline 
D-D & \textbf{\textcolor{red}{$4.31\times10^{-3}$}} & \textbf{\textcolor{red}{$154.7$}}\tabularnewline
\hline 
N-N & \textbf{\textcolor{red}{$4.7\times10^{-3}$}} & \textbf{\textcolor{red}{$294.9$}}\tabularnewline
\hline 
\end{tabular}
\par\end{centering}
\centering{}\caption{Errors and timing for solving the SPDE in equation \ref{eq:SPDE}
with filtered noise using the FIP method, for D-D and N-N boundary
conditions. The simulation parameters are $\Delta x=0.04$,\textcolor{blue}{{}
}$\Delta t=1\times10^{-3}$, and $2\times10^{4}$ samples. The errors
are dominated by the space step error. \label{tab:numerics_SPDE-DD-NN}}
\end{table}

\subsection{(1+2)-dimensional stochastic heat equation with filtered noise \label{subsec:(1+2)-d-heat-noise}}

To demonstrate the viability of the Fourier algorithm for an SPDE
with higher spatial dimension, we now extend the (1+1)-dimensional
stochastic heat equation from subsection \ref{subsec:(1+1)-d-heta_filtered_noise},
to include an additional spatial dimension, using different boundary
conditions. The (1+2)-dimensional stochastic heat equation is given
by

\begin{align}
\frac{\partial u\left(t,x,y\right)}{\partial t} & =\frac{1}{2}\left(\frac{\partial^{2}}{\partial x^{2}}+\frac{\partial^{2}}{\partial y^{2}}\right)u\left(t,x,y\right)+\eta\left(t,x,y\right).\label{eq:2d-spde}
\end{align}
We assume the noise term $\eta(t,x,y)$ is delta correlated in time,
with a noise correlation including a cutoff $\kappa$ so that:
\begin{align}
\langle\eta(t,x,y)\eta(t',x',y')\rangle & =\delta(t-t')\left(\frac{\kappa}{2}\right)^{2}\nonumber \\
 & \times\exp(-\kappa(|x-x'|+|y-y'|)).
\end{align}

We impose D-D boundary conditions on the $x$-spatial dimension and
N-N boundary conditions on the $y$-spatial dimension. The solution
$u(t,x,y)$ is expanded as

\begin{align}
u(t,x,y) & =\sum_{n\ge1,m\geq0}u_{nm}\left(t\right)\sin\left(n\pi x/L\right)\cos\left(m\pi y/L\right)\,.
\end{align}
The observable 
\begin{align}
J(t) & \equiv\intop_{0}^{L}\intop_{0}^{L}\langle u(t,x,y)^{2}\rangle\,dx\,dy\,,\label{eq:J(t)-2d}
\end{align}
assuming $u_{nm}(0)=0$, has an exact solution given by
\begin{align}
J(t) & =\sum_{n>0}I_{nn}\left[\frac{1-e^{-2K_{n0}t}}{2K_{n0}}\right]\left[1-(\frac{1-e^{-\kappa L}}{\kappa L})\right]\nonumber \\
 & +\sum_{n>0,m>0}\left[\frac{1-e^{-2K_{nm}t}}{2K_{nm}}\right]I_{nn}J_{mm}\,,
\end{align}
where 
\begin{align}
I_{nn} & =\frac{1}{1+\left(n\pi/\kappa L\right)^{2}}+\frac{2\left(n\pi/\kappa L\right)^{2}(1-e^{-\kappa L}(-1)^{n})}{\left(\kappa L\right)\left(1+\left(n\pi/\kappa L\right)^{2}\right)^{2}}\nonumber \\
J_{mm} & =\frac{1}{1+\left(m\pi/\kappa L\right)^{2}}-\frac{2(1-e^{-\kappa L}(-1)^{m})}{\left(\kappa L\right)\left(1+\left(m\pi/\kappa L\right)^{2}\right)^{2}}\,.
\end{align}

We ran numerical simulations using the FIP method, choosing the time
to range from $0$ to $0.5$, while the spatial interval was from
$0$ to $L=4$ with $\kappa=10$. The analytical solution was summed
numerically up to $n=10^{3}$ and $m=10^{3}$ to guarantee convergence.
Numerical results are tabulated in Table \ref{tab:numerics_SPDE-2D-DD:NN}.

\begin{figure}[H]
\centering{}\includegraphics[width=0.75\columnwidth]{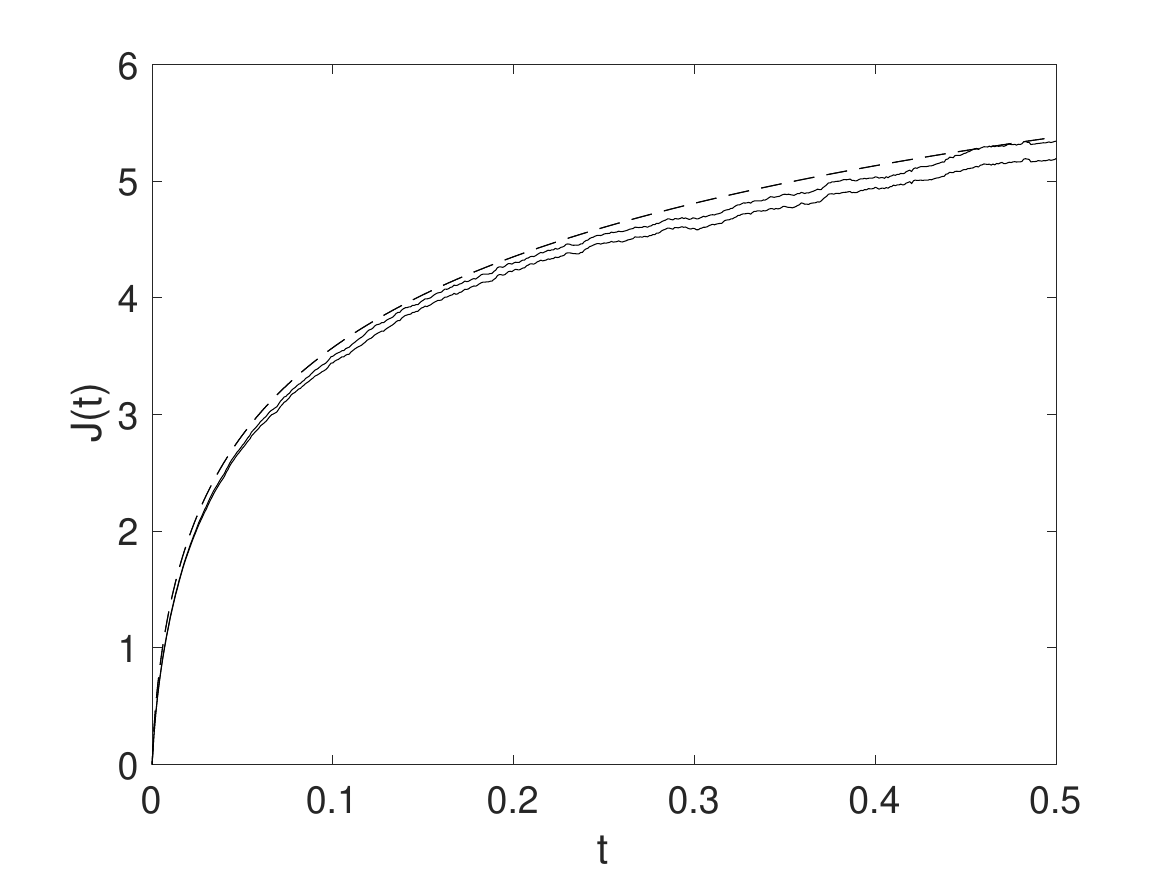}\caption{Plot of the integrated mean square solution $J(t)$ of the stochastic
heat equation with filtered noise, computed using the Fourier interaction
picture (FIP) method for the D-D boundary conditions on the $x$-spatial
dimension and N-N boundary conditions on the $y$-spatial dimension
with $\Delta x=0.05$, $\Delta y=0.05$, $\Delta t=1\times10^{-3}$
and $400$ samples. The solid lines indicate upper and lower sampling
error bounds, while the dashed line indicates the analytical solution
in Eq. (\ref{eq:J(t)-2d}) with the summations up to $n=10^{3}$ and
$m=10^{3}$.\label{fig:SPDE-2d}}
\end{figure}

\begin{table}
\begin{centering}
\begin{tabular}{|c|c|c|}
\hline 
Boundary & FIP error & Tc (s)\tabularnewline
\hline 
\hline 
D-D (x): N-N (y) & \textcolor{red}{$2.20\times10^{-2}$} & \textcolor{red}{$164.2$}\tabularnewline
\hline 
\end{tabular}
\par\end{centering}
\centering{}\caption{Errors and timing for solving the SPDE in equation \ref{eq:SPDE}
with filtered noise using the Fourier interaction picture (FIP) method.
The simulation parameters used a spatial resolution of $\Delta x=0.05$,
$\Delta y=0.05$, a time step size of $\Delta t=1\times10^{-3}$,
and a total of $400$ samples. The sampling error and time-step error
for the FIP method for these parameters are of the order of $10^{-3}$
and $10^{-2}$, respectively. \label{tab:numerics_SPDE-2D-DD:NN}}
\end{table}

\section{Conclusion\label{sec:Conclusion}}

We introduced two fast Fourier spectral methods for solving parabolic
PDEs and SPDEs with non-periodic boundaries and a uniform spatial
grid. One method evaluated derivatives using fast trigonometric transforms
(FSD), the other used an interaction picture (FIP). Both can be applied
to time-varying Dirichlet, Robin, and Neumann boundary conditions
for arbitrary spatial dimensions and field components, provided suitable
patch functions are available. We investigated a number of practical
examples. Grids were mostly chosen so that the space-step error was
dominant, in order to analyze the resulting spatial discretization
errors. Results were compared with public domain implementations of
Tau and Galerkin polynomial spectral methods.

In cases with rapidly varying spatial solutions, such as the Peregrine
and breather examples, the errors obtained using the Fourier methods
are lower than those obtained from either polynomial method. For problems
with smooth solutions relative to the spatial grid, the Fourier methods
yield spatial errors comparable to or greater than the Tau method,
due to its better scaling with spatial step-size, but less than the
Galerkin algorithm.

The FIP method is faster than the Tau and Galerkin methods for all
linear and nonlinear cases studied. The FSD method is adaptable to
a wider range of problems, and is also faster than the Tau and Galerkin
methods in almost all cases.

A scaling investigation was carried out in a nonlinear case, over
spatial grid sizes from $20$ to $640$ grid points, using very small
time-steps. Both Fourier methods gave an initially quadratic error
scaling, $\propto\Delta x^{2}$, for small $\Delta t$, while the
Tau method gave more rapid error reductions. Below the von Neumann
stability limit of $\Delta x\le\sqrt{\Delta t}$, the FIP method gave
lower errors, but at less than a quadratic rate, presumably needing
a momentum-space filter \citep{drummond1983central}. The FSD method
became unstable, while the Tau errors saturated. 

Timing results are due to a combination of algorithm, hardware and
software. We found little timing change for small grids. With larger
spatial grid size grew, there was a sub-linear timing increase for
the Fourier methods. The Tau method was the slowest, and the speed
ratio increased with grid size, with nearly an order of magnitude
speed improvement for Fourier methods at the largest size. 

Complexity analysis allows one to understand how errors can be minimized
for a given total resource use. In the Peregrine problem, using both
midpoint and RK4 time-stepping, the FIP algorithm gave lower errors
per resource use than the Tau method, although this algorithm appears
to have a slightly higher asymptotic complexity exponent.

For the linear heat equation and the stochastic heat equation, the
improvements in accuracy and speed are very substantial. The FIP method
yielded lower errors for the same resources in all cases compared
to the Tau polynomial method. We showed how the 1D stochastic heat
equation can be efficiently combined with finite noise correlation
lengths using Fourier methods, and showed how these can be used for
complex stochastic problems in higher dimensions with zero boundary
conditions, like the 2D stochastic heat equation with filtered noise. 

In summary, either Fourier or polynomial spectral methods can be utilized
for non-periodic boundary conditions, each having different areas
of suitability. There are very large efficiency improvements in linear
cases for interaction picture methods. Fourier methods are faster
than the polynomial spectral methods tested here, and have lower errors
when there is rapid spatial variation relative to the spatial grid
size. They are especially useful for stochastic partial differential
equations, spatially correlated noise, and cases with rapid spatial
variation.

\subsection*{Acknowledgements:}

We gratefully acknowledge a grant from NTT Phi Laboratories. This
publication was made possible through the support of Grant 62843 from
the John Templeton Foundation. The opinions expressed in this publication
are those of the author(s) and do not necessarily reflect the views
of the John Templeton Foundation.

\section*{Appendix}

\section*{A: Numerical approximation using spectral method}

Here we give a more detailed explanation of the method used for a
boundary value problem of the form
\begin{equation}
u_{t}=u_{xx}\label{eq:Appendix, heat equation}
\end{equation}
with initial and boundary conditions
\begin{align}
u(x,0) & =u_{0}(x)\nonumber \\
u(0,t) & =U(u(0,t))\nonumber \\
u_{x}(L,t) & =N(u(L,t)).\label{eq:Initial_boundary_conditions}
\end{align}
We first transform the function $u$ such that the boundary conditions
are homogeneous:
\begin{equation}
v=u-U(u(0,t))-xN(u(L,t))
\end{equation}

This yields a new boundary value problem, 
\begin{equation}
v_{t}=v_{xx}-\theta(x,t)
\end{equation}
with homogenous boundary conditions of the form
\begin{align}
v(0,t) & =0\nonumber \\
v_{x}(L,t) & =0,
\end{align}
with $\theta(x,t)=\dot{u}(0,t)\frac{\partial}{\partial u}U(u(0,t))+x\dot{u}(L,t)\frac{\partial}{\partial u}N(u(L,t))$.
Next we perform the discrete Fourier sine transform on $N_{s}$spatial
points,
\begin{equation}
v(x,t)=\sum_{n=1}^{N_{S}}\upsilon_{n}(t)sin(k_{n}x_{m})
\end{equation}
where $k_{n}$ is chosen such that $cos(k_{n}x)=0$. In order to obtain
a valid inverse transform, some restrictions apply to the spatial
grid $\{x_{m}\}$ on which we evaluate the discrete version of the
transform. We then obtain the derivatives
\begin{align}
v_{t}(x,t) & =\sum_{n=1}^{N_{S}}\dot{\upsilon}_{n}(t)sin(k_{n}x_{m})\nonumber \\
v_{xx}(x,t) & =-\sum_{n=1}^{N_{S}}k_{n}^{2}v_{n}(t)sin(k_{n}x_{m})
\end{align}

Taking the discrete sine transform (DST) of $v_{xx}(x,t)$ leads to
\begin{equation}
DST\left(v_{xx}(x,t)\right)=-\sum_{n=1}^{N_{S}}\sum_{l=1}^{N_{S}}k_{n}^{2}v_{n}sin(k_{n}x_{m})sin(k_{l}x_{m})
\end{equation}

If we then define the spectral derivative matrix as
\begin{equation}
D_{ln}=k_{n}^{2}\sum_{l=1}^{N_{S}}sin(k_{n}x_{m})sin(k_{l}x_{m})
\end{equation}

the derivative expression becomes
\begin{equation}
DST\left(v_{xx}(x,t)\right)=-\sum_{n=1}^{N_{S}}D_{ln}v_{n}
\end{equation}

With an appropriate choice of lattice, the derivative matrix becomes
simply diagonal. In particular, if we let $x_{m}=\left(m-\frac{1}{2}\right)\Delta x$,
with $\Delta x=L/N_{S}$, upon simplifying the trigonometric expression
for $D_{ln}$, we see that for all values of $n$ and $m,$ the only
contributing term to the sum is $n=l$, and so the derivative matrix
is diagonal, with entries
\begin{equation}
D_{nm}=k_{n}^{2}\delta_{nm}
\end{equation}

Care must be taken in implementing this propagator, since several
conventions for the discrete sine transform exist with different normalizations,
which may change the prefactor of the matrix.

\section*{B: Equivalence to Galerkin method}

We show in this Appendix that there is a relationship between the
methods we use here and the Galerkin\textcolor{black}{{} }\citep{Galerkin1968}\textcolor{black}{.}
approach to numerically solving PDEs. We consider a simple boundary
value problem of the form given in Eq (\ref{eq:Appendix, heat equation})
and \ref{eq:Initial_boundary_conditions}. That is, we seek a function
$u\in X$ that solves this problem where $X$ is the space of functions
satisfying the DN type boundary conditions.

We define a new space for test functions, $V=\left\{ v|v\in H^{1}(0,L),v(0)=0\right\} $.
Here, $H^{1}(0,L)$ denotes the first-order Sobolev space on $(0,L)$.
We multiply the equation $u_{t}=u_{xx}$ by an arbitrary function
$w\in V$. Integrating by parts leads to
\begin{align}
\intop_{0}^{L}u_{t}w(x)dx & =\intop_{0}^{L}u_{xx}w(x)dx\nonumber \\
 & =\left[u_{x}w(x)\right]_{0}^{L}-\intop_{0}^{L}u_{x}w^{'}(x)dx.\label{eq:Galerkin_int_parts}
\end{align}

Inserting the Neumann boundary condition at $x=L$ and using the definition
that $v(0)=0$ leads to
\begin{equation}
\intop_{0}^{L}u_{t}w(x)dx=N(u(L,t))w(L)-\intop_{0}^{L}u_{x}w^{'}(x)dx.
\end{equation}

We next apply the Galerkin approximation, in which the field $u$
is expanded as
\begin{align}
u_{h} & =\sum_{j\ge0}v_{j}(t)\varphi_{j}(x).
\end{align}

The term $\varphi_{0}(x)$ is chosen to satisfy the boundary condition
at $x=0$, and $v_{0}(t)=v_{0}$ is time-independent. Substituting
this approximation into Eq (\ref{eq:Galerkin_int_parts}), and choosing
the test function as $w=\varphi_{i}(x)\in B$, we obtain
\begin{align}
\sum_{j}\intop\dot{v_{j}}(t)\varphi_{j}(x)\varphi_{i}(x)dx & =N\left(\sum_{j\ge0}^ {}v_{j}(t)\varphi_{j}(x)\right)\varphi_{i}(L)\nonumber \\
 & -\sum_{j}\intop v_{j}(t)\varphi_{j}^{'}(x)\varphi_{i}^{'}(x)dx
\end{align}

Since we are dealing with a 1D problem, we may simplify this by choosing
$\varphi_{0}(x)=1$, such that $u_{0}=U$. Then one obtains a simplified
form of
\begin{align}
\sum_{j}\intop\dot{v_{j}}(t)\varphi_{j}(x)\varphi_{i}(x)dx & =N\left(U+\sum_{j\ge1}v_{j}(t)\varphi_{j}(x)\right)\varphi_{i}(L)\nonumber \\
 & -\sum_{j\ge1}\intop v_{j}(t)\varphi_{j}^{'}(x)\varphi_{i}^{'}(x)dx
\end{align}

Next we define, for $i,j>0$:
\begin{equation}
A_{ij}=\intop\varphi_{j}(x)\varphi_{i}(x)dx,
\end{equation}

\begin{equation}
B_{ij}=\intop\varphi_{j}^{'}(x)\varphi_{i}^{'}(x)dx,
\end{equation}

and
\begin{equation}
g_{i}(\mathbf{v}(t))=N\left(U+\sum_{j\ge1}v_{j}(t)\varphi_{j}(L)\right)\varphi_{i}(L).
\end{equation}

We then obtain a linear system of equations for $\dot{v}_{j}(t)$
in the form
\begin{equation}
A\mathbf{\dot{v}}(t)+B\mathbf{v}(t)=g(\mathbf{v}(t))
\end{equation}

and isolating $\dot{\mathbf{v}}(t)$ one obtains
\begin{equation}
\mathbf{\dot{v}}(t)=A^{-1}g(\mathbf{v}(t))-D\mathbf{v}(t)
\end{equation}

where $D=A^{-1}B$. In the spectral method, one may choose a basis
of half-sines, augmented with one basis function to satisfy the boundary
condition at $x=L$, where $B=\left\{ sin(k_{n}x),x|n=1.....\right\} $,
with $k_{n}=\left(n-\frac{1}{2}\right)\pi$ such that it satisfies
$cos(k_{n}x)=0$. This is essentially a Fourier sine transform, where
an extra basis function is chosen to satisfy the inhomogeneous boundaries,
while the sine functions satisfy homogeneous boundaries. The choice
of basis leads to
\begin{equation}
A_{ij}=\intop sin(k_{j}x)sin(k_{i}x)dx=\frac{1}{2}\delta_{ij}
\end{equation}

and:
\begin{equation}
B_{ij}=k_{i}k_{j}\intop cos(k_{j}x)cos(k_{i}x)dx=\frac{k_{i}k_{j}}{2}\delta_{ij}=\frac{k_{i}^{2}}{2}\delta_{ij}.
\end{equation}

Finally, one obtains
\begin{equation}
D_{ij}=k_{i}^{2}\delta_{ij}.
\end{equation}

The matrix $D$ is exactly equal to the spectral derivative matrix.
The spectral method is therefore a case of the Galerkin approximation
method, restricted to homogeneous boundary conditions. We generalize
this to inhomogeneous boundaries through a patch technique, in order
to maintain the efficient use of fast Fourier transforms.

\section*{C: DST and DCT transforms}

The types of discrete sine and cosine transforms used depend on the
boundary conditions, and are defined below. Here $N$ is the number
of spatial points including points at boundaries, and $N_{T}$ is
the points in the DCT/DST transform. This can differ, since any zero
boundary values are not required in the transforms. We also define
$N_{FT}$ as the Fourier transform \citep{Frigo1998,Frigo2005} logical
size, giving the number of spatial points used in a discrete Fourier
transform implementation of the DCT/DST transform.

As explained in the main text, in (\ref{subsec:Boundary-types}),
there can be multiple boundary types for each dimension and field
component. In boundary case ``$\beta$'', the forward and inverse
transforms $\mathcal{F}^{\beta}$ and $\mathcal{I}^{\beta}$ are such
that in all cases, for an arbitrary field component $u$ sampled at
discrete points, with a discrete transform $\tilde{u}$, one has that:
\begin{align}
\tilde{u}_{n} & =\sum_{j}\mathcal{F}_{nj}^{\beta}u_{j}\nonumber \\
u_{j} & =\sum_{m}\mathcal{I}_{jn}^{\beta}\tilde{u}_{n}.\label{eq:Forward-backward-transforms}
\end{align}

In all cases the grid spacing in the $i-th$ dimension in position
space and momentum space is:
\begin{align}
\Delta x_{i} & =\frac{X_{i}^{b}-X_{i}^{a}}{N-1}\nonumber \\
\Delta k_{i} & =\frac{\pi}{X_{i}^{b}-X_{i}^{a}}.
\end{align}

The boundaries are labelled by the $\beta$ index describing the type.
We use standard DCT/DST transform names \citep{Strang1999,Pueschel2003,Frigo2005}
to identify the forward and inverse transform pairs, on the basis
that each transform pair has to satisfy the specified zero boundary
conditions on either the function or the derivative, as described
in the main text. In all cases the corresponding number of discrete
Fourier transform points is $N_{FT}=2(N-1)$.

\subsection*{1: D-D case: DST-I}

Let $u(0)=u_{1}=0$, and $u(R)=u_{N}=0$. The discrete representation
of $u$ is:

\paragraph{Forward transform: DST-I}

\begin{align}
\tilde{u}_{n} & =\sqrt{\frac{2}{N-1}}\sum_{j=2}^{N-1}u_{j}\sin\left(\pi\frac{\left(j-1\right)\left(n-1\right)}{N-1}\right).
\end{align}

\paragraph{Inverse transform: DST-I
\begin{align}
u_{j} & =\sqrt{\frac{2}{N-1}}\sum_{n=2}^{N-1}\tilde{u}_{n}\sin\left(\pi\frac{\left(j-1\right)\left(n-1\right)}{N-1}\right).
\end{align}
}

The forward transform does \textbf{not} require the values at the
end-points which are set to zero in this case. This is implicit in
the sine expansion, since $\text{\ensuremath{\sin}(n\ensuremath{\pi)=0}.}$
Second derivatives are proportional to $\left(n-1\right)^{2}$.

\subsection*{2: D-N case: DST-II/III}

Let $u(0)=0$, and $u'(R)=0$. The discrete representation of $u$
is:

\paragraph{Forward transform: DST-III}

\begin{align}
\tilde{u}_{n} & =\sqrt{\frac{2}{N-1}}\Biggl((-1)^{\left(n-1\right)}u_{N}/2\nonumber \\
 & \left.\ +\sum_{j=2}^{N-1}u_{j}\sin\left[\frac{\pi}{N-1}\left(n-\frac{1}{2}\right)(j-1)\right]\right).
\end{align}

\paragraph{Inverse transform: DST-II}

\begin{align}
u_{j} & =\sqrt{\frac{2}{N-1}}\left(\sum_{n=1}^{N-1}\tilde{u}_{n}\left(t\right)\sin\left[\frac{\pi}{N-1}\left(n-\frac{1}{2}\right)j\right]\right).
\end{align}

The forward transform does \textbf{not} require the value at $n=1$,
which is zero in this case. This is implicit in the sine expansion,
since $\text{\ensuremath{\sin}(n\ensuremath{\pi)=0}.}$ This transform
implies a boundary condition that is odd around $n=1$, and even around
$n=N$. Second derivatives are proportional to $\left(n-1/2\right)^{2}$.

\subsection*{3: N-D case: DCT-II/III}

Take $u'(0)=u(R)=0$. The discrete representation of $u$ is:

\paragraph{Forward transform: DCT-III}

\begin{align}
\tilde{u}_{n} & =\sqrt{\frac{2}{N-1}}\Biggl(u_{1}/2\nonumber \\
 & \left.\ +\sum_{j=2}^{N-1}u_{j}\cos\left[\frac{\pi}{N-1}\left(n-\frac{1}{2}\right)(j-1)\right]\right).
\end{align}

\paragraph{Inverse transform: DCT-II}

\begin{align}
u_{j} & =\sqrt{\frac{2}{N-1}}\sum_{n=1}^{N}\tilde{u}_{n}\cos\left[\frac{\pi}{N-1}\left(n-\frac{1}{2}\right)\left(j-1\right)\right].
\end{align}

The forward transform does not require the value at $n=N$, which
is zero in this case. This transform implies a boundary condition
that is even around $j=1$, and odd around $j=N$. Second derivatives
are proportional to $\left(n-1/2\right)^{2}$.

\subsection*{4: N-N case: DCT-I}

Let $u'(0)=0$, and $u'(R)=0$. The discrete representation of $u$
is:

\paragraph{Forward transform: DCT-I}

\begin{align}
\tilde{u}_{n} & =\sqrt{\frac{2}{N-1}}\Biggl(\frac{1}{2}u_{1}+\frac{1}{2}(-1)^{n-1}u_{N}\nonumber \\
 & \left.\ +\sum_{j=2}^{N-1}u_{j}\cos\left(\pi\frac{\left(j-1\right)\left(n-1\right)}{N-1}\right)\right).
\end{align}

\paragraph{Inverse transform: DCT-I}

\begin{align}
u_{j} & =\sqrt{\frac{2}{N-1}}\Biggl(\frac{1}{2}\tilde{u}_{1}+\frac{1}{2}(-1)^{j-1}\tilde{u}_{N}\nonumber \\
 & \left.\ +\sum_{n=2}^{N-1}\tilde{u}_{n}\cos\left(\pi\frac{\left(j-1\right)\left(n-1\right)}{N-1}\right)\right).
\end{align}

The forward transform requires the values at the end-points, which
are not zero in this case. It is equal (up to a factor) to a discrete
Fourier transform of $2\left(N-1\right)$ real numbers with even symmetry
about $j=1$ and $j=N$. As a result, the equivalent discrete derivatives
at both the end-points are zero. Second derivatives are proportional
to $\left(n-1\right)^{2}$.

\section*{D: Error analysis}

Here, we define the RMS relative comparison error used in this work
for all comparisons. All examples have analytical solutions $a(t_{i},x_{j})$,
which are compared with the computed solutions $c(t_{i},x_{j})$,
by taking their absolute differences \textcolor{black}{$d\left(t_{i},x_{j}\right)=|c\left(t_{i},x_{j}\right)-a\left(t_{i},x_{j}\right)|$,
}for all time and position points. The RMS relative comparison error
$\epsilon$ is then \textcolor{black}{
\begin{equation}
\epsilon_{}=\sqrt{\frac{\sum_{i}\sum_{j}d(t_{i},x_{j})^{2}}{N_{T}N_{S}m}}\,,
\end{equation}
where $N_{T}$ is the number of time points, $N_{S}$ is the number
of spatial points, and $m$ is the maximum absolute value of the }computed
result. If the spatial grid is non-uniform, the square of each result
\textcolor{black}{$d\left(t_{i},x_{j}\right)$ should be weighted
by the space-step ($\varDelta x_{j}$) $\times$ time-step ($\varDelta t_{}$)
and then the total divided by $XT$, not the total space-time points
($N_{j}N_{T}$), giving}

\textcolor{black}{
\begin{equation}
\epsilon_{}=\sqrt{\frac{\sum_{i}\sum_{j}d(t_{i},x_{j})^{2}\Delta x_{j}}{XN_{T}m}}
\end{equation}
}

\textcolor{black}{where $\Delta x_{j}$ is the spatial step size, $X$ is
the range of $x$, with $N_{T}=T/\Delta t$, where $T$ is the range
of time. Since the spatial grid generated is non-uniform, the spatial
step size $\Delta x_{j}=x_{j}-x_{j-1}$ varies for each step, while
the time step remains the same for all time points.} We use the term
'error' to define the RMS relative comparison error, $\triangle t$ for
the time step size and $\triangle x$ for the spatial step size.

\section*{E: Comparison between finite difference and spectral methods\label{sec:Comparison-between-finite}}

We solve the heat equation and the nonlinear Schrödinger equation
in one space dimension with zero boundaries as a test to compare the
error results of the finite difference (FD) method with those of the
FIP and FSD methods. The finite difference (FD) method used here is
the central differencing method and implemented in xSPDE4 \citep{Drummond2023}.

In both examples, time-steps were chosen to be small enough to give
negligible time-discretization errors when using the FIP method. This
isolates the spatial discretization errors. In some cases the time-step
used is not short enough to prevent instabilities with other methods.
The errors are still reported in the tabulated results for comparison
purposes. The spatial discretization error can be reduced by using
a finer spatial grid. The tests are carried out with an iterated midpoint
(MP) method, to improve the accuracy of the finite difference comparisons.

\subsection*{1: Heat equation - zero boundary case \label{subsec:Heat-equation}}

The equation for the (1+1)-dimensional heat equation is given by

\begin{eqnarray}
\frac{\partial u}{\partial t} & = & \frac{\partial^{2}u}{\partial x^{2}}\,.
\end{eqnarray}

To compare the errors here we conduct two tests of the heat equation
using two combinations of boundary conditions: Dirichlet-Dirichlet
(D-D) and Neuman-Neuman (N-N). For simplicity, we consider the boundary
values to be set to zero for both Dirichlet and Neuman boundary conditions.

\subsubsection*{Dirichlet-Dirichlet (D-D) \label{subsec:Dirichlet-Dirichlet-(D-D)}}

We consider a time-dependent solution of the form $u(t,x)$ that satisfies
the boundary condition $u(0,0)=u(0,1)=0$ in the $x$ interval $\left[0,\pi\right]$
and $t$ interval $\left[0,1\right]$. The exact solution takes the
form
\begin{equation}
u(t,x)=\sum_{n=1}^{\infty}c_{n}\sin\left(nx\right)e^{-n^{2}t}\ ,
\end{equation}
where $c_{n}=2\int_{0}^{1}u(x,0)\sin\left(nx\right)dx$. For this
case one possible solution takes the form
\begin{equation}
u(t,x)=2\sin\left(x\right)e^{-t}+\sin\left(2x\right)e^{-4t},
\end{equation}
 Average (RMS) relative errors and computational times for the two
Fourier spectral methods and the finite difference method are shown
in \textcolor{black}{Table \ref{tab:HeatEqnDD}.}

\begin{table}[H]
\centering{}%
\begin{tabular}{|c|c|c|c|c|c|c|}
\hline 
$\Delta t$ & FSD error & Tc (s) & FIP error & Tc (s) & FD error & Tc (s)\tabularnewline
\hline 
\hline 
$0.001$ & \textcolor{red}{$\mathbf{9\times10^{-8}}$} & \textcolor{red}{$\mathbf{0.44}$} & \textbf{\textcolor{red}{$\mathbf{5\times10^{-14}}$}} & \textcolor{red}{$\mathbf{0.33}$} & $6\times10^{-4}$ & $0.27$\tabularnewline
\hline 
$0.002$ & \textcolor{red}{$\mathbf{3\times10^{-7}}$} & \textcolor{red}{$\mathbf{0.23}$} & \textbf{\textcolor{red}{$\mathbf{2\times10^{-14}}$}} & \textcolor{red}{$\mathbf{0.19}$} & $6\times10^{-4}$ & $0.14$\tabularnewline
\hline 
$0.1$ & \textcolor{red}{$\mathbf{\infty}$} & \textcolor{red}{$\mathbf{0.02}$} & \textbf{\textcolor{red}{$\mathbf{4\times10^{-16}}$}} & \textcolor{red}{$\mathbf{0.02}$} & $\mathbf{\infty}$ & $0.02$\tabularnewline
\hline 
\end{tabular}\caption{Comparison of the error and elapsed computational time (Tc) in solving
the Heat equation with D-D boundary condition using the Fourier spectral
derivative (FSD), Fourier interaction picture (FIP) and finite difference
(FD) methods for three different time step sizes $\Delta t$. The
number of spatial steps is fixed at 20, with $\Delta x=\pi/20$.\label{tab:HeatEqnDD}}
\end{table}

\textcolor{black}{In this case, the FIP solution is exact up to machine
precision. With the FIP method, the error is reduced by more than
a factor of $10^{10}$ compared to the FD method and $10^{6}$ compared
to the FSD method. We observe that for the same time-step, the FIP
method is slower compared to the FD method but faster than the FSD
method. However, there is no requirement to use extremely small $\Delta t$,
as shown in Table \ref{tab:HeatEqnDD}. The use of relatively large
time steps, such as $\Delta t=0.1$, with the FIP method would still
result in small errors of the order of machine precision ($\sim10^{-16}$)
for this particular problem. In fact from the Table \ref{tab:HeatEqnDD},
one sees that an advantage of the FIP method is that accurate results
up to machine precision can be obtained using large time steps, significantly
reducing the elapsed time.}

\subsubsection*{Neumann-Neumann (N-N)}

Next, we consider a time dependent solution of the form $u(t,x)$
that satisfies the boundary condition $u^{'}(0,0)=u^{'}(0,1)=0$ in
the $x$ interval $\left[0,\pi\right]$ and t interval $\left[0,1\right]$.
The exact solution is
\begin{equation}
u(t,x)=\sum_{n=1}^{\infty}d_{n}cos\left(nx\right)e^{-n^{2}t}\ ,
\end{equation}
where $d_{n}=2\int_{0}^{1}u(x,0)cos\left(nx\right)dx$. For this case
one possible solution takes the form

\begin{equation}
u(t,x)=2+cos\left(x\right)e^{-t}+cos\left(2x\right)e^{-4t}.
\end{equation}
Relative errors and computational times for the three methods considered
are shown in \textcolor{black}{Table \ref{tab:HeatEqnNN}.}

\begin{table}[H]
\centering{}%
\begin{tabular}{|c|c|c|c|c|c|c|}
\hline 
$\Delta t$ & FSD error & Tc (s) & FIP error & Tc (s) & FD error & Tc (s)\tabularnewline
\hline 
\hline 
$0.001$ & \textcolor{red}{$\mathbf{6.0\times10^{-8}}$} & \textcolor{red}{$\mathbf{0.47}$} & \textcolor{red}{$\mathbf{1\times10^{-14}}$} & \textcolor{red}{$\mathbf{0.36}$} & $4\times10^{-4}$ & $0.28$\tabularnewline
\hline 
$0.002$ & \textcolor{red}{$\mathbf{2.3\times10^{-7}}$} & \textcolor{red}{$\mathbf{0.24}$} & \textcolor{red}{$\mathbf{6\times10^{-15}}$} & \textcolor{red}{$\mathbf{0.18}$} & $4\times10^{-4}$ & $0.15$\tabularnewline
\hline 
$0.1$ & \textcolor{red}{$\mathbf{\infty}$} & \textcolor{red}{$\mathbf{0.02}$} & \textbf{\textcolor{red}{$\mathbf{1\times10^{-16}}$}} & \textcolor{red}{$\mathbf{0.02}$} & $\mathbf{\infty}$ & $0.03$\tabularnewline
\hline 
\end{tabular}\caption{Comparison of the error obtained and computation time $Tc$ for solving
the heat equation with N-N boundary conditions using the Fourier spectral
derivative (FSD), Fourier Interaction picture (FIP) and finite difference
(FD) methods for three different time step sizes $\Delta t$. The
number of spatial steps is fixed at 20, with $\Delta x=\pi/20$.\label{tab:HeatEqnNN}}
\end{table}

\textcolor{black}{Using the FIP method for this problem reduces the
error by more than a factor of $10^{10}$ compared to the FD method
and $10^{6}$ compared to the FSD method. The results obtained with
N-N boundary conditions are similar to those seen with using }D-D\textcolor{black}{{}
boundary conditions. Specifically, using the FIP Method leads to errors
that are close to the limits of machine precision as shown in Table
\ref{tab:HeatEqnNN}. Additionally, when the time step size is increased,
the errors are of the same order or smaller, significantly reducing
the overall runtime.}

\subsection*{2: Nonlinear Schrödinger equation - zero boundary case}

The equation for the (1+1)-dimensional nonlinear Schrödinger equation
is given by

\begin{eqnarray}
\frac{\partial u}{\partial t} & = & i\cdot\left(u\cdot\left|u\right|^{2}-\frac{u}{2}+\frac{1}{2}\frac{\partial^{2}u}{\partial x^{2}}\right)\,.
\end{eqnarray}

This has the exact, time-independent solution of $u=sech\left(x\right)$.
For a finite interval test, one can use the exact solution as an initial
condition and then displace the solution to give a zero Dirichlet
boundary condition. This is achieved by choosing $\tilde{u}\left(x,t\right)=u\left(x,t\right)-u_{0}$,
with the equation:

\begin{eqnarray}
\frac{\partial\tilde{u}}{\partial t} & = & i\cdot\left(u\cdot\left|u\right|^{2}-\frac{u}{2}+\frac{1}{2}\frac{\partial^{2}\tilde{u}}{\partial x^{2}}\right)\,.
\end{eqnarray}
For this case the solution takes the form of $\tilde{u}=sech\left(x\right)-u_{0}$.
To obtain an exactly soluble problem with a zero Neumann boundary
condition for testing purposes, one boundary is chosen at $x=0$,
which is at the center of the $sech\left(x\right)$ solution.

\begin{table}[H]
\begin{centering}
\begin{tabular}{|c|c|c|c|c|c|c|}
\hline 
Bdary & FSD error & Tc(s) & FIP error & Tc(s) & FD error & Tc(s)\tabularnewline
\hline 
\hline 
D-D & \textcolor{red}{$\mathbf{9\times10^{-5}}$} & \textbf{\textcolor{red}{$\mathbf{0.47}$}} & \textbf{\textcolor{red}{$\mathbf{9\times10^{-5}}$}} & \textbf{\textcolor{red}{$\mathbf{0.34}$}} & $3\times10^{-4}$ & 0.33\tabularnewline
\hline 
D-N & \textcolor{red}{$\mathbf{2\times10^{-5}}$} & \textbf{\textcolor{red}{$\mathbf{0.72}$}} & \textbf{\textcolor{red}{$\mathbf{2\times10^{-5}}$}} & \textbf{\textcolor{red}{$\mathbf{0.46}$}} & $9\times10^{-5}$ & 0.39\tabularnewline
\hline 
\end{tabular}
\par\end{centering}
\centering{}\caption{Comparison of the error for solving the non-linear Schrödinger equation
with D-D and D-N boundary conditions using the Fourier spectral derivative
(FSD), Fourier interaction picture (FIP) and finite difference (FD)\textcolor{black}{{}
methods} for a fixed spatial and time step size. For all methods $\Delta x=0.1$
and $\Delta t=0.001$.\label{tab:Table_NLS}}
\end{table}

Results in Table \ref{tab:Table_NLS} indicate a 25-fold error reduction
for D-D and a 4-fold error reduction for D-N boundary combinations
for the FSD and FIP methods, compared to the FD method.

\subsubsection*{3: Summary}

In summary, for the comparisons of spatial errors with Fourier and
finite difference methods, we find that the Fourier methods give much
lower errors, with slightly longer computational times per step. However,
the number of steps in time can be greatly reduced in the FIP method
for linear cases, and this gives both faster times and lower errors.\\

\bibliographystyle{unsrt}

\end{document}